\documentclass[10pt,reqno]{amsart}
\usepackage{fullpage}
\usepackage{amsfonts}
\usepackage{amssymb,amsthm,amsmath}
\usepackage{times}
\usepackage{textgreek}
\usepackage{imakeidx}
\usepackage{graphicx}
\usepackage{float}
\usepackage[title]{appendix}
\usepackage{subcaption}
\usepackage{mathrsfs}
\usepackage{appendix}
\usepackage{hyperref}

\hypersetup{
	colorlinks=true,
	linkcolor=black,
	filecolor=mangenta,
	urlcolor=cyan,
	citecolor=blue,
	pdftitle={Limit Theorems For Non-Hermitian Ensembles},
	bookmarks=true,
	pdfpagemode=FullScreen}

\usepackage[UKenglish]{datetime}
\newdateformat{UKvardate}{\THEDAY\ \monthname[\THEMONTH] \THEYEAR}

\newtheorem{theorem}{Theorem}[]
\newtheorem*{theorem*}{Theorem}

\newtheorem{lemma}[theorem]{Lemma}

\newtheorem{corollary}[theorem]{Corollary}

\newtheorem{remark}{Remark}[section]
\newtheorem*{remark*}{Remark}
\hypersetup{colorlinks=true, linkcolor=black, filecolor=black, urlcolor=black }
\DeclareMathOperator{\Tr}{Tr}
\DeclareMathOperator{\erf}{erf}

\begin{document}
\small \title{Limit Theorems For Non-Hermitian Ensembles} 
\maketitle
\begin{center}
\uppercase \author{Olivia V. Auster}\\ \text{}\\
\address{School of Mathematical Sciences \\ 
Queen Mary University of London \\ Mile End Road \\ London, E1 4NS \\ 
United Kingdom \\ \text{} \\
\small \email{o.v.auster@qmul.ac.uk\\ \text{}\\}}
\end{center}

\begin{abstract} The complex Ginibre ensemble and its generalisation, the complex induced Ginibre ensemble, have been well-studied in the field of Random Matrix Theory. In the present work, the asymptotic distribution and the independence of the extreme eigenvalue moduli are studied in the limit of large dimensions of random matrices for these non-Hermitian ensembles. They are derived with the use of Andreief’s integration formula and the known methodological approach defined for the study of the limiting distribution of the scaled spectral radius at the edge of the complex Ginibre ensemble. The limiting distribution of the scaled spectral radius and the scaled minimum modulus for the complex induced Ginibre ensemble, with a proportional rectangularity index different from zero, is the Gumbel distribution. In the limit of a large size of complex Ginibre matrices, the left and right tails of the distribution of the minimum modulus are the Rayleigh and Weibull distributions, respectively. The limiting left tail of the distribution of the minimum modulus is the same for these random matrix ensembles, with a rectangularity index of the complex induced Ginibre ensemble equal to zero. This phenomenon is also verified for the right tail of the distribution of this random minimum. The independence of the extreme moduli is formally established, at appropriate scaling, for large matrices from the complex Ginibre ensemble and the complex induced Ginibre ensemble with a proportional rectangularity index. This study extends knowledge in the field of Random Matrix Theory for random variables, like the minimum modulus of matrices from non-Hermitian ensembles, whose limiting stochastic dynamics have never been explored.
\end{abstract}

\maketitle

\section{Introduction}\label{Introduction}\text{}\\
The attractiveness of the Random Matrix Theory (RMT) lies in the possibility of using it as a means to model problems in high dimensions and perform related calculations analytically. This is mainly due to the invariance property of the probability distribution of certain random matrix ensembles. Matrix ensembles that are characterised by an invariant probability distribution are those from classical compact matrix groups studied in references \cite{Weyl1946} and \cite{Meckes2019}. 
The complex Ginibre ensemble is one of them and is a special case of the Ginibre-Girko ensemble with maximal non-Hermiticity as presented in \cite{FyodorovSommersKhoruzhenko1998} and \cite{AkemannPhillipsShifrin2009}. Like the real and the symplectic Ginibre ensembles, this ensemble was introduced by J. Ginibre \cite{Ginibre1965}. The distribution of eigenvalues of a complex Ginibre matrix is comparable to that of the distribution of the positions of charges of a two-dimensional Coulomb gas in a harmonic oscillator potential, at a specific temperature corresponding to the Dyson index $\beta =2$.\\ \\ 
The complex Ginibre ensemble is the space of $N \times N$ complex matrices $J$, whose complex entries are the $J_{ij} = x_{ij} + \text{i} y_{ij} \in \mathbb{C}$, that are independent and identically distributed.  Each of them follows a standard complex Gaussian distribution $\mathcal{N}_{\mathbb{C}}(0,1)$ with probability density function 
\begin{equation}
	P(J_{ij})= \frac{1}{\pi}e^{-\vert J_{ij} \vert^2}
\end{equation}
where $(x_{ij})_{1\leq i,j \leq N} \in \mathbb{R}^{N \times N}$ and $(y_{ij})_{1\leq i,j \leq N} \in \mathbb{R}^{N \times N}$. The real and imaginary parts of the entries denoted $x_{ij}$ and $y_{ij}$ are independent random variables, each following a real Gaussian distribution $\mathcal{N}_{\mathbb{R}}(0, \frac{1}{2})$.\\ 
\\ 
The space of $N \times N $ complex matrices from the complex Ginibre ensemble is endowed with a probability measure here denoted $\mu(J)$, where $d\mu(J) = P(J) \vert D(J)\vert$ and the joint probability density function of the entries $J_{ij}$ is 
\begin{align}
	  P(J) = \frac{1}{\pi^{N^2}}\exp\left(-\Tr\left( J^{\dagger}J\right)\right)
\end{align}
The term $\vert D(J) \vert = \frac{1}{2^{N^2}}\bigwedge_{k=1}^{N}\bigwedge_{j=1}^{N} \vert dJ_{kj}^{*}\wedge dJ_{kj}\vert$ is called the Cartesian volume element as presented in \cite{KhoruzhenkoSommers2015}.\\
\\
The joint probability density function of the eigenvalues of an $N \times N$ complex Ginibre matrix $J$ is expressed as
\begin{equation}\label{EigenvaluesJointProbabilityDensityFunction}
	P_{N}(z) = \frac{1}{\pi^{N}C_{N}}\exp\left(-\sum_{k=1}^{N}\vert z_k\vert^2 \right) \vert \Delta (z)\vert^2 
\end{equation}
where $z = \left\lbrace z_1, \dots, z_{N} \right\rbrace$ is the set of the complex eigenvalues of the random matrix $J$ and $\Delta(z)$ is the Vandermonde determinant $\Delta (z) = \prod_{1 \leq i < j \leq N} \left( z_i - z_j \right)$, with $\vert \Delta (z)\vert^2 = \Delta (z) \overline{\Delta (z)}$. The term $C_N$ is the normalisation factor defined as
\begin{equation}
	C_N = \frac{1}{\pi^N}\int_{\mathbb{C}^N}\exp \left( -\sum_{k=1}^{N}\vert z_k\vert^2 \right)  \prod_{1 \leq i < j \leq N}\vert z_i - z_j \vert^2\prod_{k=1}^{N}d^{2}z_{k} = \prod_{k=1}^{N}\Gamma(k+1)
\end{equation}
\\
A topic which has been widely investigated in RMT is the limiting distribution of the eigenvalues of $N \times N$ random matrices as $N$ goes to infinity. Studies of the limiting distribution of the largest and smallest eigenvalues for Gaussian ensembles are presented in \cite{TracyWidom1994, TracyWidom1996, TracyWidom2009}. The distributions of eigenvalues (and their spacings) have been studied in \cite{AkemannPhillipsShifrin2009} with the computation of gap probabilities with respect to radial ordering for non-Hermitian random matrices and their chiral counterparts. The statistical properties of the spectral radius of matrices from non-Hermitian ensembles have also been studied in the literature, as it is for the complex Ginibre ensemble \cite{Rider2003, Rider2004} and the real Ginibre ensemble \cite{RiderSinclair2014}.\\
\\
The present paper is structured as follows. Section \ref{Methods} details the methodology used. The latter encompasses the application of Andreief's integration formula \cite{Forrester2018} and the framework established by B. Rider in \cite{Rider2003}. The results are presented in Section \ref{Results}. Section~\ref{LimitTheoremComplexGinibre} is devoted to the analysis of the distribution of the modulus of the extreme eigenvalues of matrices from the complex Ginibre ensemble. The limiting left and right tails of the distribution of the minimum modulus are studied for large matrices in Section~\ref{LimitingDistributionRmin}. Pursuing the analysis for this random matrix ensemble in the scaling limit $\sqrt{N}$ as $N$ goes to infinity, the independence of the scaled minimum modulus with respect to the scaled spectral radius is formally established in Section~\ref{IndependenceScaledRminSacledRmaxComplexGinibre}. The eigenvalues of matrices from the complex Ginibre ensemble form a determinantal point process. The joint probability density of their respective modulus has been investigated in \cite{Kostlan1992}. This result is extended in \cite{HoughKrishnapurPeresVirag2006}.
Employing the method introduced in \cite{Rider2003}, from which limit theorems are derived at the edge of the spectrum of the complex and symplectic Ginibre ensembles, results presenting the limiting distribution of the scaled spectral radius and minimum modulus for the complex induced Ginibre ensemble \cite{FischmannBruzdaKhoruzhenkoSommersZyczkowski2012}, are derived in Section~\ref{LimitTheoremComplexInducedGinibre}. As an ideal modelling framework for the expansion of the research related to non-Hermitian matrices, the complex induced Ginibre ensemble is a special case of the Feinberg-Zee ensemble and a generalisation of the complex Ginibre ensemble. The joint probability density function of the entries of matrices $A$ from the Feinberg-Zee ensemble is 
\begin{equation}
	P_{FZ}(A) \propto  \exp \left( -\Tr V(A^{\dagger}A) \right)
\end{equation}
The induced Ginibre ensembles correspond to the Feinberg-Zee ensemble with potential $V(y) = -\frac{\beta}{2} ( y - L\log(y))$, where $L$ is a non-negative parameter called the rectangularity index, and $\beta = 1$ or $\beta = 2$ for the real or complex induced Ginibre ensembles, respectively.\\ \\ Let $G$ denote an $N \times N$ complex induced Ginibre matrix. The joint probability density function of its entries is 
\begin{equation}
	P_{indG}^{(\beta)}(G) = \Gamma^{(\beta)} \left(\det \left(G^{\dagger}G \right)\right)^{\frac{\beta}{2}L}\exp \left( -\frac{\beta}{2}\Tr G^{\dagger}G \right) 
\end{equation}
where $\beta = 2$ and the determinant $\left(\det \left(G^{\dagger}G\right) \right)^{a} = \prod_{k=1}^{N} \vert \lambda_{i} \vert^{2a}$, $a \geq 0$. The factor $\Gamma^{(\beta)} $ is the partition function
\begin{equation}
	\Gamma^{(\beta)} = \pi^{-\frac{\beta}{2}N^{2}}\left( \frac{\beta}{2}\right)^{\frac{1}{2}N(N-L)} \prod_{j=1}^{N}\frac{\Gamma(\frac{\beta}{2}j)}{\Gamma(\frac{\beta}{2}(j + L)}
\end{equation}
\\
In the limit of large matrix dimensions, the eigenvalues are spread across an annulus in the complex plane, which is distinguished from the phenomenon identified for the complex Ginibre ensemble, thus referenced as the circular law and presented in the work by V. L. Girko \cite{Girko1985} and the paper by T. Tao, V. Vu and M. Krishnapur \cite{TaoVuKrishnapur2010}. The joint probability density function of the eigenvalues of an $N \times N$ complex induced Ginibre matrix $G$ is 
\begin{equation}\label{ProbabilityDensityFunctionComplexInducedGinibre}
	P_{N}(\lambda_{1}, \cdots, \lambda_{N}) = \frac{1}{Z_{N}} \prod_{j<k}^{N}\vert  \lambda_{k} - \lambda_{j}\vert^{2} \prod_{j=1}^{N}\vert \lambda_{j} \vert^{2L}\exp\left(- \sum_{j=1}^{N}\vert \lambda_{j} \vert^{2} \right)
\end{equation}
with partition function 
\begin{equation}
	Z_{N} 
		 = \int  \cdots \int \prod_{j<k}^{N} \vert  \lambda_{k} - \lambda_{j}\vert^{2} \prod_{j=1}^{N}\vert \lambda_{j} \vert^{2L}\exp\left(-\sum_{j=1}^{N}\vert \lambda_{j} \vert^{2} \right)\prod_{j=1}^{N}d^{2}\lambda_{j}
		 = N! \pi^{N} \prod_{k=1}^{N} \Gamma(k + L)
\end{equation}
It is understood from the definition of the joint probability density function of the eigenvalues that there is a repulsion of the eigenvalues from the origin of the complex plane due to the term $\left(\det(G^{\dagger}G)\right)^{L} = \prod_{j=1}^{N}\vert \lambda_{j} \vert^{2L}$ appearing in the formula. More precisely, the distribution of the eigenvalues is supported by a ring centred at the origin of the complex plane with an outer circle of radius $r_{out} = \sqrt{L+N}$ and an inner circle of radius $r_{in} = \sqrt{L}$. This phenomenon refers to the "Single Ring Theorem" stating that, in the limit of large matrix dimensions, the empirical distribution of the eigenvalues has for support a single ring \cite{FeinbergZee1997} \cite{GuionnetKrishnapurZeitouni2011}. For a small rectangularity index $L$ (close to zero), the mean eigenvalue density converges to the mean eigenvalue density for the complex Ginibre ensemble.
\\
A study is performed in Section~\ref{LimitTheoremComplexInducedGinibre} for the complex induced Ginibre ensemble from the work of B. Rider \cite{Rider2003}. In this reference, the author details an analytical method to determine the limiting distribution of the scaled spectral radius of matrices from the complex and symplectic Ginibre ensembles. A similar analysis is undertaken in Section~\ref{LimitTheoremsOuterInnerEdge} for the scaled spectral radius and scaled minimum modulus of an $N \times N$ complex induced Ginibre matrix. The study is conducted at the outer and inner edges of the ring defined with a rectangularity index  assumed to be proportional to $N$, i.e., $L = \alpha N$, $\alpha > 0$. An exact fit between the empirical distribution and analytical formulation of the distribution of the extreme moduli is acknowledged numerically. Additionally, the independence of the (scaled) spectral radius and (scaled) minimum modulus is studied in the limit as $N$ goes to infinity for these two non-Hermitian ensembles in Section~\ref{IndependenceExtremeModuliComplexInducedGinibre}.
\\
The limiting tails of the distribution of the eigenvalue moduli are also of interest in the field. The similitude of the limiting right tail of the minimum modulus distribution between the two non-Hermitian ensembles under scrutiny is explored in Section \ref{ComparisonsLimitingDistribution}. The same analysis is performed for the left tail of the distributions of this extreme modulus. A discussion and a list of conclusions are set out in Section \ref{Conclusions}.\\

\section{Methods}\label{Methods}\text{}\\
Let $A$ denote an $N \times N$ complex Ginibre matrix. Let $P\left( r_{min}^{(N)}(A) \geq a \right)$ denote the survival distribution function of the minimum modulus $r_{min}^{(N)}(A)$ and $P\left( r_{max}^{(N)}(A) \leq a \right)$ the cumulative distribution function of the spectral radius $r_{max}^{(N)}(A)$. These distribution functions are derived analytically in the present paper with the use of Andreief's integration formula presented in \cite{Forrester2018}. The joint probability density function of the eigenvalues of the matrix $A$, denoted $P_N(z_1, \cdots, z_N)$, is defined as in Section \ref{Introduction} with Equation~(\ref{EigenvaluesJointProbabilityDensityFunction}).
The probability that the minimum modulus $r_{min}^{(N)}(A)$ is greater than the radius $a \in \mathbb{R}^{+*}$, is 
\begin{align*}
	P\left(r_{min}^{(N)}(A) \geq a \right) 
	& = \int_{\vert z_1\vert \geq a} \cdots \int_{\vert z_N\vert \geq a}P_N(z_1, ..., z_N)\prod_{k =1}^N d^2z_k\\
	& = \frac{1}{Z_N}\int_{\vert z_1 \vert >a}\cdots \int_{\vert z_N \vert >a} 
\det\left[ z_k^{N-j} \right]_{j,k=1}^N \det\left[ \bar{z}_k^{N-j} \right]_{j,k=1}^N\prod_{k = 1}^N dm(z_k)
\end{align*}
where $dm(z_k) = e^{-\vert z_k\vert^2}d^2z_k$ and $Z_{N} = N! \pi^{N}\prod_{k=0}^{N-1} k!$ is the partition function.\\
\\
Then, performing a change of variables in polar coordinates gives, 
\begin{equation*}
P\left(r_{min}^{(N)}(A) \geq a \right) 
				= \frac{1}{Z_N} N!\prod_{k=0}^{N-1} 2\pi\int_{a}^{+ \infty} e^{- r^2}r^{2k} rdr
				= \prod_{k=0}^{N-1}\frac{1}{\Gamma(k+1)}\int_{a^2}^{+ \infty} e^{- t}t^k dt
\end{equation*}\\
This probability corresponds to the probability that no eigenvalue lies inside the disk of radius $a$. It is analytically formulated as
\begin{equation}\label{SurvivalPrRmin}
	P\left(r_{min}^{(N)}(A)\geq a \right) = \prod_{k=0}^{N-1}\frac{\Gamma(k+1, a^{2})}{\Gamma(k+1)}
\end{equation}\\
This result is initially stated in \cite{GrobeHaakeSommers1988} without proof and without scaling applied to the minimum modulus. Asymptotic expansions are studied in \cite{Forrester1992}. Equation~(\ref{SurvivalPrRmin}) is also presented in \cite{AkemannPhillipsShifrin2009}, in which Gram's formula \cite{Mehta2004} is applied to get the probability that no eigenvalue lies in the disk of radius $a$ centred at the origin of the complex plane as a finite product of regularised upper incomplete Gamma functions. 
It is straightforward to derive the probability density function $p_{r_{min}^{(N)}(A)}$ of the minimum modulus $r_{min}^{(N)}(A)$ from Equation (\ref{SurvivalPrRmin}). More precisely, it is convenient to take the logarithm of the probability $P(r_{min}^{(N)}(A) \geq r)$, $r \in \mathbb{R}^{+*}$, and find its first derivative with respect to $r$. Let the survival distribution function $P(r_{min}^{(N)}(A) \geq r)$ denote the function $H(r)$ of the radius $r$. The corresponding probability density function of the minimum modulus $r_{min}^{(N)}(A)$ is $p_{r_{min}^{(N)}(A)}(r) = -\frac{\text{d}H(r)}{\text{d} r}$.\\ \\
Furthermore, 
\begin{align*}
	\frac{\text{d} H(r)}{\text{d} r} 
		& = H(r)\frac{\text{d}}{\text{d} r}\left(\log\left[H(r)\right]\right) 
		  = \prod_{k=0}^{N-1}\frac{\Gamma(k+1, r^{2})}{\Gamma(k+1)} \left[ \sum_{k=0}^{N-1}\frac{\partial}{\partial r}  \left[ \log\left( \Gamma(k+1, r^{2})\right) \right]\right]
\end{align*}
with,
\begin{align*}
	\frac{\partial\log\left[ \Gamma(k+1, r^{2})\right]}{\partial r}  
		& = \frac{1}{\Gamma(k+1, r^{2})} \frac{\partial \Gamma(k+1, r^{2})}{\partial r}
\end{align*}
It is known that $\frac{\partial \Gamma(s,x)}{\partial x} = -x^{s-1}e^{-x}$.\\
\\ 
Thus,
\begin{align*}
	\frac{\text{d} H(r)}{\text{d} r} 
	& = \prod_{k=0}^{N-1}\frac{\Gamma(k+1, r^{2})}{\Gamma(k+1)} \sum_{j=0}^{N-1}\left[ \frac{1}{\Gamma(j+1, r^{2})}  \frac{\partial \Gamma(j+1, r^{2})}{\partial r}  \right]
\end{align*}
where
\begin{equation*}
	\frac{\partial \Gamma(j+1, r^{2})}{\partial r} = -2r^{2j+1}e^{-r^{2}} 
\end{equation*}\\
Finally, the probability density function of the minimum modulus $r_{min}^{(N)}(A)$ is, \\
\begin{align}
	p_{r_{min}^{(N)}(A)}\left(r\right) 
	& = 2re^{-r^{2}}\prod_{k=0}^{N-1}\frac{\Gamma(k+1, r^{2})}{\Gamma(k+1)}\sum_{j=0}^{N-1}\left[  \frac{r^{2j}}{\Gamma(j+1, r^{2})}  \right]
\end{align}
\\
Using a similar approach as for the case of no scaling, the survival distribution function of the scaled minimum modulus $\frac{r_{min}^{(N)}(A)}{\sqrt{N}} $ for an $N\times N$ complex Ginibre matrix $A$ has the same analytical expression as for $r_{min}^{(N)}(A)$ replacing $r$ with $\sqrt{N}r$.
The corresponding probability density function is the probability density function of the minimum modulus $r_{min}^{(N)}(A)$ multiply by $\sqrt{N}$ at the point $\sqrt{N}r$, i.e.,  
\begin{equation}
	p_{\frac{r_{min}^{(N)}(A)}{\sqrt{N}}}\left(r\right)
	 = \sqrt{N} p_{r_{min}^{(N)}(A)}\left(\sqrt{N} r\right)
\end{equation}
\\ \\
The distribution function the spectral radius, denoted $r_{max}^{(N)}(A)$, is explored in the following. 
\\
\begin{lemma}\label{CumulativePrRmaxLemma} Let $A$ denote an $N \times N$ complex Ginibre matrix. Let $r_{max}^{(N)}(A)$ denote the spectral radius.\\ \\ Then,
\begin{equation}\label{CumulativePrRmax}
	P\left(r_{max}^{(N)}(A) \leq r \right) = \prod_{k=0}^{N-1}\frac{\gamma(k+1, r^2)}{\Gamma(k+1)}
\end{equation}
where $\Gamma(k)$ is the Gamma function and $\gamma(k,r)$ is the lower incomplete Gamma function.
\end{lemma} 
\begin{proof} The cumulative distribution function of the spectral radius $r_{max}^{(N)}(A)$ for the complex Ginibre ensemble is easily derived with the use of Andreief's integration formula. More precisely,
\begin{align*}
	P\left(r_{max}^{(N)}(A) \leq r \right) 
	& = \int_{\vert z_1\vert \leq r} \cdots \int_{\vert z_N\vert \leq r}P_N(z_1, \cdots, z_N)\prod_{k =1}^N d^2z_k 
	  = \prod_{k=0}^{N-1}\frac{1}{\Gamma(k+1)}\int_{0}^{r^2}e^{-t}t^{k}dt
\end{align*}\\
\end{proof}

\begin{corollary} The joint probability law of the moduli $ \vert z_{1} \vert, ..., \vert z_{N} \vert$ of an $N \times N$ complex Ginibre matrix is the joint distribution of independent random variables $\gamma_{k}$ each following a Gamma-Rayleigh distribution $GR(\alpha_{k}, \delta_{k})$, with $\alpha_{k} = k$ and $\delta_{k} = 1$, i.e., we have the equality in distribution
\begin{equation}
		\left\lbrace \vert z_{1} \vert, ..., \vert z_{N} \vert \right\rbrace  \overset{d}{=} \left\lbrace \gamma_{1}, ..., \gamma_{N} \right\rbrace 
\end{equation}
where $\gamma_{k} ~ GR(k, 1)$.\\
\end{corollary}
\begin{proof}
	\begin{equation*}\label{ProductGammaRayleighDistribution}
	 	P(r_{1}^{(N)}(A) \leq r, r_{2}^{(N)}(A) \leq r, ...,  r_{N}^{(N)}(A) \leq r) 
			= \prod_{k=1}^{N} P(\gamma_{k} < r)
\end{equation*}
where the random variables $\gamma_{k}$ are independent. Each random variable $\gamma_{k}$ follows a Gamma-Rayleigh distribution $GR(\alpha_{k}, \delta_{k})$,\\ with parameters $\alpha_{k} = k$ and $\delta_{k} = 1$. 
\end{proof}\text{}\\
The Gamma-Rayleigh distribution is derived in the work by E. Akarawak, I. Adeleke and R. Okafor \cite{AkarawakAdelekeOkafor2018}. \\ \\
\text{}\\
The corresponding probability density function of the largest modulus $r_{max}^{(N)}(A)$ is derived from Equation (\ref{CumulativePrRmax}). Taking the logarithm of the cumulative distribution function $P(r_{max}^{(N)}(A) \leq r)$ and then the first derivative with respect to $r$ of $\log\left(P(r_{max}^{(N)}(A) \leq r) \right)$, the probability density function of the spectral radius $r_{max}^{(N)}(A)$ is,\\
\begin{equation}
	p_{r_{max}^{(N)}(A)}(r) 
		= 2re^{-r^{2}}\prod_{k=0}^{N-1}\frac{\gamma(k+1, r^{2})}{\Gamma(k+1)}\sum_{j=0}^{N-1}\left[  \frac{r^{2j}}{\gamma(j+1, r^{2})}  \right]
\end{equation} \\
\\
A similar method of derivation, as the one applied for the scaled minimum modulus, leads to the formulation of the probability density function of the scaled spectral radius for the complex Ginibre ensemble as 
\begin{equation}
	p_{\frac{r_{max}^{(N)}(A)}{\sqrt{N}}}\left(r\right) 
	= \sqrt{N} p_{r_{max}^{(N)}(A)}\left( \sqrt{N} r\right)
\end{equation}
\\
Asymptotic results are then derived from these analytical expressions in the limit as the size of matrices goes to infinity. They are presented in Section ~\ref{Results}.\\ 
\\
The limiting distribution of the scaled spectral radius $\frac{r_{max}^{(N)}(A)}{\sqrt{N}}$ of an $N \times N$ complex Ginibre matrix is established in the research paper by B. Rider \cite{Rider2003}. His framework aims to derive the logarithm of the cumulative distribution function of the scaled spectral radius $\frac{r_{max}^{(N)}(A)}{\sqrt{N}}$ for the complex Ginibre ensemble, in the limit as $N$ goes to infinity. The limiting distribution of the extreme moduli $r_{min}^{(N)}(G)$ and $r_{max}^{(N)}(G)$, as well as their independence, for an $N \times N$ matrix $G$ from the complex induced Ginibre ensemble, are derived at appropriate scaling with his framework in Section~\ref{Results}. The main steps are detailed as follows. First, analytical formulae of the distribution function of these random variables are derived as a partial product of the distribution function of a weighted sum of independent and identically distributed random variables, each following a standard exponential distribution. As an illustration, B. Rider \cite{Rider2003} considered the scaled spectral radius $\frac{r_{max}^{(N)}(A)}{\sqrt{N}}$ and established that 
\begin{align}
	P\left(\frac{r_{max}^{(N)}(A)}{\sqrt{N}} \leq r \right) 
	  = \prod_{k=0}^{N-1}P\left( \frac{1}{N}\sum_{j=1}^{N-k}Z_{j} \leq r^{2} \right)
\end{align}
where the random variable $Z^{(k+1)} = \sum_{j=1}^{k+1}Z_{j}$ follows a Gamma distribution with shape parameter $k+1$ and rate parameter $1$. The random variables $Z_{j}$, with $j \in \lbrace 1, ..., k+1 \rbrace$, are independent and identically distributed. Each of the $Z_{j}$ follows a standard exponential distribution. Then, he derived lower and upper bounds for this analytical formula in the limit as $N$ goes to infinity. The limit, as $N$ goes to infinity, of the logarithm of this distribution function is then established in his paper by applying the classical Edgeworth expansion. Then, an approximation from this expression considering the leading order sum as the Riemann integral is derived. Finally, he chose an appropriate definition of a function $f_{N}$ such that the limiting distribution of the considered random variable, here the scaled spectral radius $\frac{r_{max}^{(N)}(A)}{\sqrt{N}}$, is the identified statistical extreme distribution, here the Gumbel distribution.\\ \\
\newpage
\section{Results}\label{Results}\text{}
\subsection{Limit theorems for the complex Ginibre ensemble}\label{LimitTheoremComplexGinibre}\text{}\\ 
\\
The limiting distributions of the extreme moduli of matrices from non-Hermitian ensembles are of interest in the field. The limiting distribution of the spectral radius is studied in references \cite{Rider2003, Rider2004, RiderSinclair2014}. The limiting left and right tails of the distribution of the minimum modulus as well as the independence of the extreme moduli have not been investigated in the literature and are the themes of the next Sections \ref{LimitingDistributionRmin}, \ref{IndependenceScaledRminSacledRmaxComplexGinibre}, \ref{LimitTheoremsOuterInnerEdge} and \ref{IndependenceExtremeModuliComplexInducedGinibre} for the complex Ginibre ensemble and its matrix ensemble generalisation, the complex induced Ginibre ensemble. Limiting stochastic expressions of these extreme random variables are established for the complex induced Ginibre ensemble in Section~\ref{LimitTheoremsOuterInnerEdge}. \\
 
\subsubsection{Limiting left and right tails of the distribution of the minimum modulus for the complex Ginibre ensemble as $N$ goes to infinity}\label{LimitingDistributionRmin}\text{}\\
\begin{theorem}\label{RminComplexGinibreRayleighTheorem}
Let $A$ denote an $N \times N$ matrix from the complex Ginibre ensemble. The left tail of the distribution of minimum modulus $r_{min}^{(N)}(A)$ is the Rayleigh distribution with parameter $\sigma = \frac{1}{\sqrt{2}}$, as $N$ goes to infinity. More precisely, for $0 < r \ll 1$, 
\begin{equation}
	\lim_{N\longrightarrow + \infty}P(r_{min}^{(N)}(A) < r) =  1 -  e^{-r^{2}}\left( 1 - O\left(r^{4} \right) \right)
\end{equation}
\end{theorem}
\begin{proof} The statement of Theorem \ref{RminComplexGinibreRayleighTheorem} is derived from the formulation of the probability $P(r_{min}^{(N)}(A) \geq r)$ which corresponds to the $N$-th partial product presented in Equation~(\ref{SurvivalPrRmin}),
\begin{align*}
	\prod_{k=0}^{N-1}\frac{\Gamma(k+1, r^{2})}{\Gamma(k+1)} 
	& =  e^{-r^{2}}\prod_{k=1}^{N-1} e^{-r^{2}}e_{k}(r^{2}) 
	  = e^{-r^{2}}H^{(N)}(r,0)
\end{align*}\\
where $e_{n}(x) = \sum_{k=0}^{n}\frac{x^{k}}{k!}$ defines the $n$-th partial sum of the exponential function. $H^{(N)}(r,0)$ is the conditional probability that given one eigenvalue lies at the origin of the complex plane all the others are found outside the disk centred at zero with radius $r$.\footnote{The derivation of the conditional probability $H^{(N)}(r,0)$ is presented in the appendix.} 
\\
\\
Thus,
\begin{align}
	\lim_{N \rightarrow +\infty} P(r_{min}^{(N)}(A) \geq r) 
		= e^{-r^{2}}\lim_{N \rightarrow +\infty}H^{(N)}(r,0)
\end{align}\\
The conditional probability $H^{(N)}(r,0)$ is a result presented in \cite{GrobeHaakeSommers1988} and \cite{KhoruzhenkoSommers2015}, where, in the limit as $N$ goes to infinity and $r$ is small (or $r$ close to zero)
\begin{equation}\label{LimitConditionalProbability}
	\lim_{N \rightarrow +\infty} H^{(N)}(r,0) =  1 - \left[\frac{r^{4}}{2} + \frac{r^{6}}{6} + \frac{r^{8}}{24} + O(r^{10}) \right] = 1 - O\left( r^{4} \right)
\end{equation}
\\
The limit presented in Equation (\ref{LimitConditionalProbability}) is derived as follows. The lower incomplete gamma function has series expansion, for $x$ around zero, 
\begin{equation}
	\gamma(a,x) = x^{a}\sum_{n=0}^{+\infty}(-1)^{n}\frac{x^{n}}{n! (a+n)} 
\end{equation}
Thus,
\begin{equation*}
	\frac{\gamma(k+1, r^{2})}{\Gamma(k+1)} = \frac{r^{2(k+1)}}{(k+1)\Gamma(k+1)} + \frac{ r^{2(k+1)}}{\Gamma(k+1)}\sum_{n=1}^{+\infty}(-1)^{n}\frac{r^{2n}}{n! (k+n+1)}
\end{equation*}\\
Then, for small $r$, i.e., $ 0 < r \ll 1$,
\begin{equation*}
	 \lim_{N \rightarrow +\infty} H^{(N)}(r,0) = \lim_{N \rightarrow +\infty}\prod_{k=1}^{N-1}\left[ 1 - \frac{r^{2(k+1)}}{(k+1)\Gamma(k+1)} + O\left( r^{2(k+2)} \right) \right]
\end{equation*}
\\
Thus,\\
\begin{align*}
	\lim_{N \rightarrow +\infty}\prod_{k=1}^{N-1} \frac{\Gamma(k+1, r^{2})}{\Gamma(k+1)}
	& = \exp\left[  \sum_{k=1}^{+\infty} \left[  - \sum_{\gamma =1}^{+\infty}\frac{1}{\gamma}\left( \frac{r^{2(k+1)}}{(k+1)!} + O\left( r^{2(k+2)} \right) \right)^{\gamma}\right] \right]\\
	& = \sum_{j=0}^{+\infty} \frac{(-1)^{j}}{j!}\left[  \sum_{k=1}^{+\infty} \left( \frac{r^{2(k+1)}}{(k+1)!} + O\left(  r^{2(k+2)}\right) \right) \right]^{j}\\
\end{align*}
which implies for small $r$ (or $r$ close to zero),
\begin{equation*}
	\lim _{N \rightarrow +\infty}  H^{(N)}(r,0) = 1 - O \left( r^{4}  \right)
\end{equation*}
Consequently, for $0 < r \ll 1$,
\begin{equation*}
	\lim_{N \rightarrow +\infty} P\left( r_{min}^{(N)}(A) < r \right)  = 1 -  e^{-r^{2}}\left( 1 - O\left(r^{4} \right) \right)
\end{equation*}\\ \\
This corresponds to the cumulative distribution of the Rayleigh distribution with parameter $\sigma = \frac{1}{\sqrt{2}}$. \\
\end{proof}

\begin{theorem}\label{RminComplexGinibreWeibullTheorem}
Let $A$ denote an $N \times N$ matrix from the complex Ginibre ensemble. The right tail of the distribution of the minimum modulus $r_{min}^{(N)}(A)$ is the Weibull distribution with shape parameter $\kappa = 4$ and scale parameter $\lambda  = \kappa^{\frac{1}{\kappa}}$ as $N$ goes to infinity, i.e., for large $r$,
\begin{equation}	
	\lim_{N \longrightarrow +\infty}P\left( r_{min}^{(N)}(A) < r \right) = 1 -  e^{-\frac{r^{4}}{4}\left( 1 + O\left( \frac{1}{ r^{2}} \right) \right)} \text{}
\end{equation}
\end{theorem}
\begin{proof} This result is established from Equation $9$ of \cite{GrobeHaakeSommers1988}.\\
\end{proof}\text{}\\

\subsubsection{Independence of the scaled spectral radius and the scaled minimum modulus for the complex Ginibre ensemble}\label{IndependenceScaledRminSacledRmaxComplexGinibre}\text{}\\
\begin{theorem} Let A denote an $N \times N$ complex Ginibre matrix. The scaled spectral radius $R_{N} = \frac{r_{max}^{(N)}(A)}{\sqrt{N}}$ and the scaled minimum  modulus $r_{N} = \frac{r_{min}^{(N)}(A)}{\sqrt{N}}$ are independent random variables, in the limit as $N$ goes to infinity.
\end{theorem}
\begin{proof} Applying Andreief's integration formula \cite{Forrester2018},
\begin{align*}
	P\left(r_{N} \geq r \text{  and  } R_{N} \leq R \right)
		& = \prod_{k=0}^{N-1}\frac{\gamma(k+1, NR^{2})}{\Gamma(k+1)}  \prod_{k=0}^{N-1}\left[ 1 -  
			\frac{\gamma(k+1, Nr^{2})}{\gamma(k+1, NR^{2})}\right]
\end{align*}
\\
Furthermore, taking the limit as $N$ goes to infinity
\begin{align*}
	& \lim_{N \rightarrow +\infty} P\left(r_{N} \geq r \text{  and  } R_{N} \leq R \right)
	  = \lim_{N \rightarrow +\infty}  \prod_{k=0}^{N-1}\frac{\gamma(k+1, NR^{2})}{\Gamma(k+1)}	
	\lim_{N \rightarrow +\infty} \prod_{k=0}^{N-1}\left[ 1 - 
	\frac{\gamma(k+1, Nr^{2})}{\gamma(k+1, NR^{2})}\right]
\end{align*}
\\
It is known that
\begin{equation*}
	\lim_{N \rightarrow +\infty}P\left( R_{N} \leq R  \right) = \lim_{N \rightarrow +\infty}  \prod_{k=0}^{N-1}\frac{\gamma \left(k+1, N R^{2}\right)}{\Gamma(k+1)} 
\end{equation*}
\\
Furthermore,
\begin{align*}
	\prod_{k=0}^{N-1}\left[ 1 -  \frac{\gamma(k+1, Nr^{2})}{\gamma(k+1, NR^{2})}\right] 
		& = \prod_{k=0}^{N-1}\left[ 1 -  \frac{\gamma(k+1, Nr^{2})}{\Gamma(k+1)} \right]  \prod_{k=0}^{N-1} \left[  1 + E_{k,N}(R,r)  \right]\\
\end{align*}\text{}\\
where
\begin{equation*}
	E_{k,N}(R,r)  = \frac{\gamma(k+1, Nr^{2})}{\Gamma^{2}(k+1)} O\left(\Gamma\left(k+1, NR^{2}\right)\right) +  \frac{\gamma^{2}(k+1, Nr^{2})}{\Gamma^{3}(k+1)} O\left(\Gamma\left(k+1, NR^{2}\right)\right)
\end{equation*}\\ 
\\
\textit{\textbf{Lower bound.}} For any $r$ and $R\in \mathbb{R}^{+*}$,
\begin{align*}
	\prod_{k=1}^{N} \left( 1 +  E_{k,N}(R,r)  \right)	
	& > \prod_{k=1}^{N} \left( 1 + \frac{\gamma(k+ 1, Nr^{2})O\left(\Gamma\left(k+1, NR^{2}\right)\right)}{N^{2}\Gamma(k+1)^{2}} \right)  
	  = 1 + O\left(\frac{1}{N} \right) 
\end{align*}\\ \\ 
\\
\textit{\textbf{Upper bound.}}  Let $\varepsilon_{k,N}(R) = O\left(\Gamma\left(k+1, NR^{2}\right)\right)$. Let $x= NR^{2}$. This implies, for any $R\in \mathbb{R}^{+*}$, that \\
\begin{align*}
	\prod_{k=1}^{N} \left( 1 + E_{k,N}(R,r)   \right)
	& = \prod_{k=1}^{N} \left( 1 +  \frac{\gamma(k+1, Nr^{2})\varepsilon_{k,N}(R)}{\Gamma(k+1)^{2}}  + O\left( \frac{\gamma^{2}(k+1, Nr^{2})\varepsilon_{k,N}(R)}{(\Gamma(k+1))^{3}}\right)\right)\\
	& = \prod_{k=1}^{N} \left( 1  + O\left( \frac{\gamma(k+1, Nr^{2})\varepsilon_{k,N}(R)}{(\Gamma(k+1))^{2}}\right)\right)\\
	& < \prod_{k=1}^{N} \left( 1  + O\left( \frac{\varepsilon_{k,N}(R)}{\Gamma(k+1)}\right)\right)\\
	& = \exp \left[ \sum_{k=1}^{N}\log \left( 1  + O\left( \frac{\varepsilon_{k,N}(R)}{\Gamma(k+1)}\right)\right) \right]\\
	& = \exp \left[ \sum_{k=1}^{N} O\left( \frac{\varepsilon_{k,N}(R)}{\Gamma(k+1)}\right) \right]
	  = \exp \left[ O\left(e^{-x}\sum_{j=0}^{k}\frac{ x^{j+1}}{j!} \right) \right]
\end{align*}\\
\\
where $x$ goes to infinity, as $N$ goes to infinity. This implies that $e^{-x}\sum_{j=0}^{k}\frac{ x^{j+1}}{j!}$ converge to zero.\\ 
\\
\\
Applying the squeeze theorem, 
\begin{equation*}
	\lim_{N \rightarrow + \infty} \prod_{k=1}^{N} \left( 1 + E_{k,N}(R,r)  \right) = 1
\end{equation*}\\
\\
Thus,
\begin{equation*}
	\lim_{N \rightarrow +\infty} \prod_{k=0}^{N-1}\left[ 1 -  \frac{\gamma(k+1, Nr^{2})}{\gamma(k+1, NR^{2})}\right] 
	= \lim_{N \rightarrow +\infty}  P\left(r_{N} \geq r \right)
\end{equation*}
\\
Consequently,
\begin{align*}
	\lim_{N \rightarrow +\infty} P\left(r_{N} \geq r \text{  and  } R_{N} \leq R \right) 
		& = \lim_{N \rightarrow +\infty}P\left( r_{N} \geq r \right)\lim_{N \rightarrow +\infty}P\left( R_{N} \leq R  \right) 
\end{align*}\\
The scaled spectral radius $R_{N}$ and the scaled minimum modulus $r_{N}$ for the complex Ginibre ensemble are then independent random variables in the limit as $N$ goes to infinity.\\
\end{proof}\text{}\\

\subsection{Limit theorems for the complex induced Ginibre ensemble}\label{LimitTheoremComplexInducedGinibre}\text{}\\ \\
The limiting distribution of the spectral radius of non-Hermitian ensembles as well as its precise localisation near the edge of the unit disk is investigated in the work by B. Rider \cite{Rider2003}. He dedicated his studies to the complex and symplectic Ginibre ensembles introduced by J. Ginibre \cite{Ginibre1965}. These statistical ensembles share similar features such as the universality conjecture known as the circular law \cite{TaoVuKrishnapur2010}. The following results are derived for the complex induced Ginibre ensemble using the methodological approach exposed in B. Rider's work \cite{Rider2003}.\\ \\ Let $G$ denote an $N \times N$ random matrix from the complex induced Ginibre ensemble with rectangularity index $L$ as defined in \cite{FischmannBruzdaKhoruzhenkoSommersZyczkowski2012}. The scaled spectral radius $R_{N} = \frac{r_{max}^{(N)}(G)}{\sqrt{L+N}}$ evolves near the outer radius $r_{out}$ while the scaled minimum modulus $r_{N} = \frac{r_{min}^{(N)}(G)}{\sqrt{L}}$ is a random variable fluctuating around the inner radius $r_{in}$. 
\\
\subsubsection{Limit theorems at the outer and inner edges}\label{LimitTheoremsOuterInnerEdge} 
\begin{figure}[H]	
	\includegraphics[scale=.45]{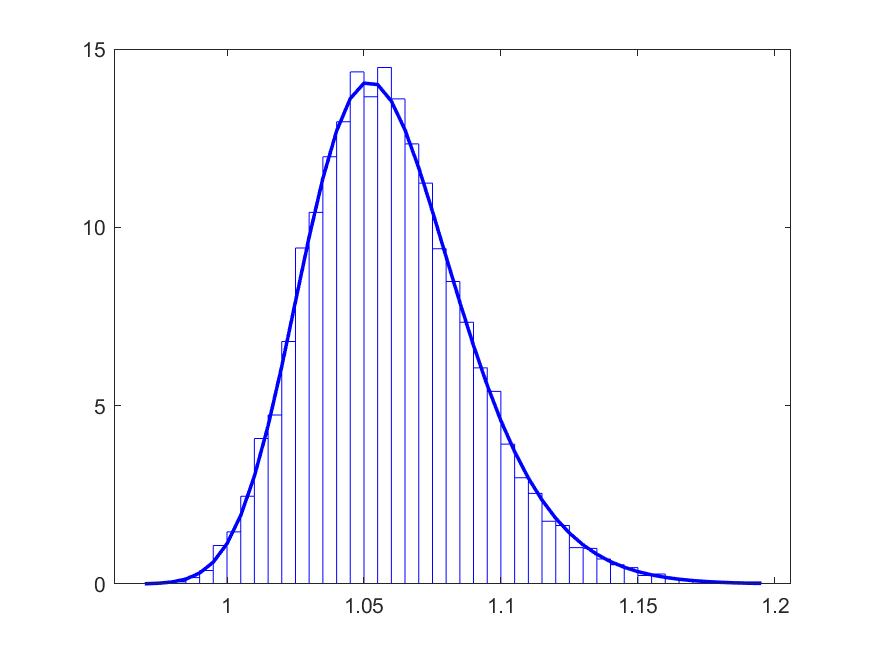}
	\hfill
	\includegraphics[scale=.45]{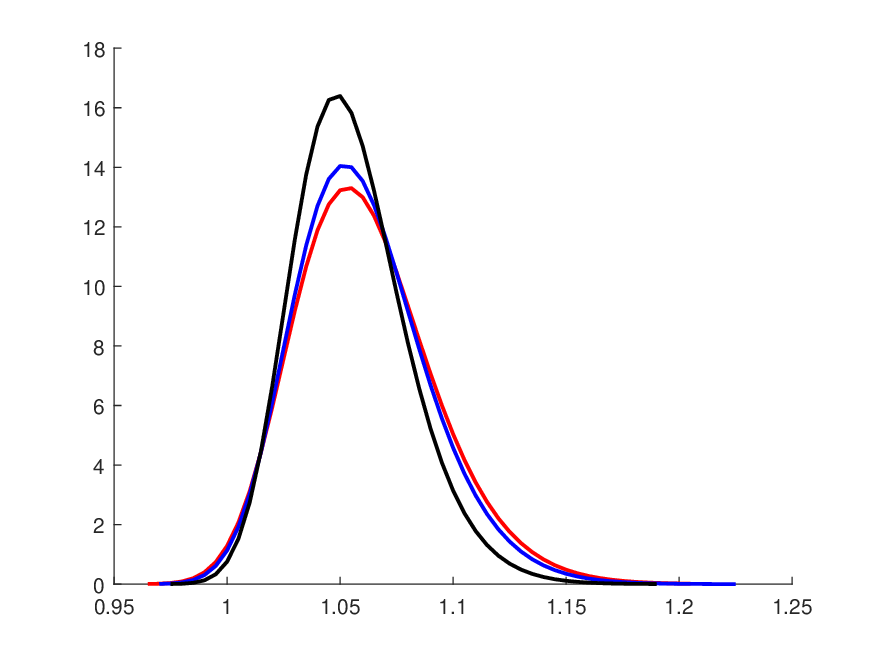}
	\caption{\footnotesize \textbf{(left panel)} Empirical probability density function (histogram) of the scaled spectral radius $R_{N}$ for $K=10^{4}$ generated $N \times N$ matrices from the complex induced Ginibre ensemble with $N =90$ and a proportional rectangularity index $L = \alpha N$, $\alpha = \frac{1}{9}$. The corresponding exact (analytical) probability density function $p_{R_{N}}$ (solid curve). \textbf{(right panel)} The analytical (exact) probability density function of the scaled spectral radius $R_{N}$ for $K=10^{4}$ generated $N \times N$ matrices from the complex induced Ginibre ensemble with $N = 90$. The rectangularity index $L$ is proportional to $N$ such that $L = \alpha N$, $\alpha > 0$. The results are presented for different values of $\alpha = \left\lbrace \frac{1}{90} , \frac{1}{9}, \frac{4}{9} \right\rbrace$.}
	\label{ComplexInducedGinibreScaledSpectralRadiusRmaxProportionalL}
\end{figure}

\begin{figure}[H]
	\begin{center}
			\includegraphics[scale=0.5]{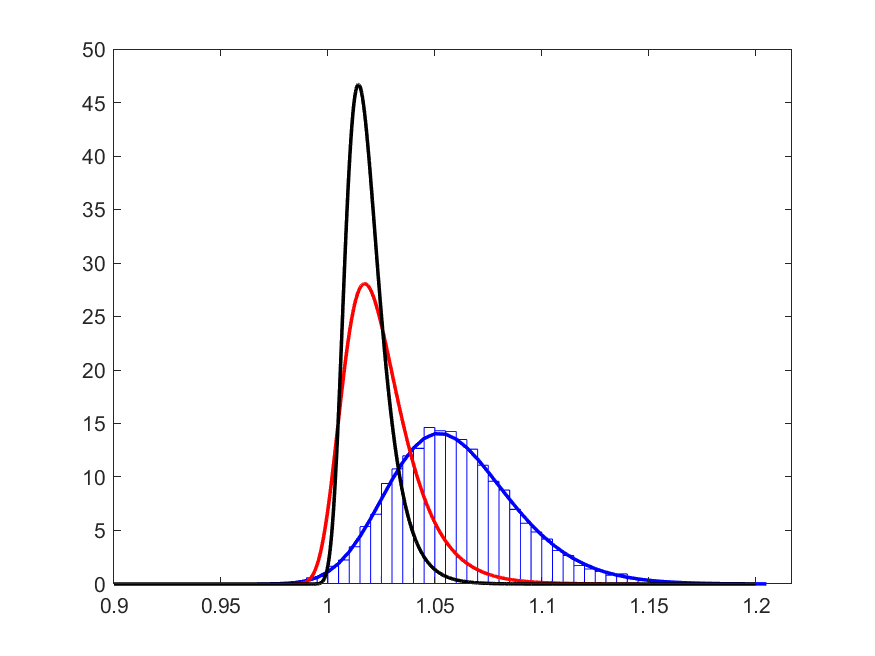}			
	\end{center}	
	\caption{\footnotesize Empirical probability density function (histogram) of the scaled spectral radius $R_{N}$ for $K=10^{4}$ generated $N \times N$ matrices from the complex induced Ginibre ensemble with $N =90$ and a proportional rectangularity index $L = \alpha N$, $\alpha = \frac{1}{9}$. The exact (analytical) probability density function $p_{R_{N}}$ (blue curve). Limiting probability density functions of the scaled spectral radius $R_{N}$ represented by the red and black curves for large $N = 10^{3} $ and $N = 2 \times 10^{3}$, respectively.}	\label{ComplexInducedGinibreLargeNRmaxDistribution}
\end{figure}

\begin{figure}[H]
	\begin{center}
		\includegraphics[scale=.45]{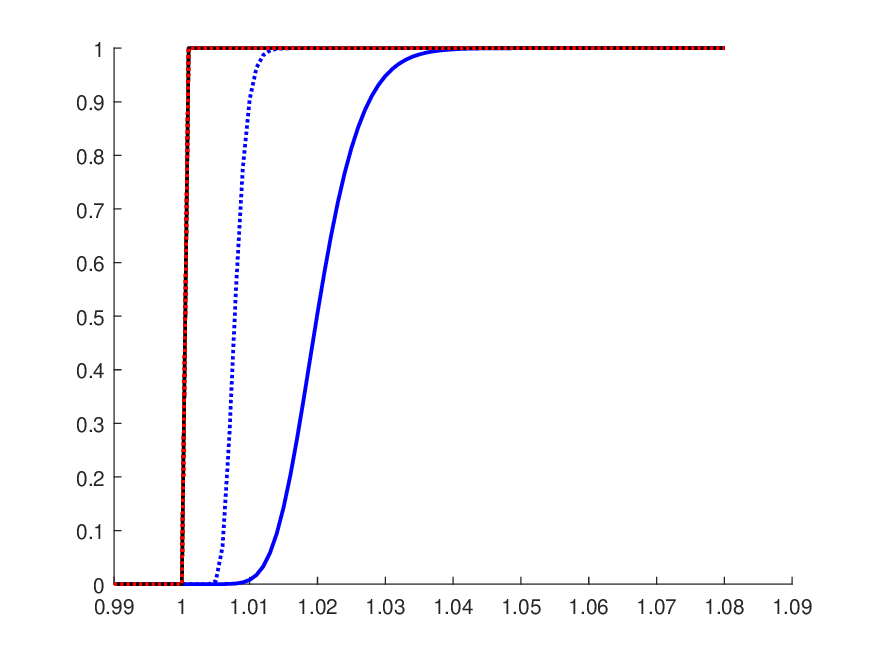}
		\hfill
		\includegraphics[scale=.45]{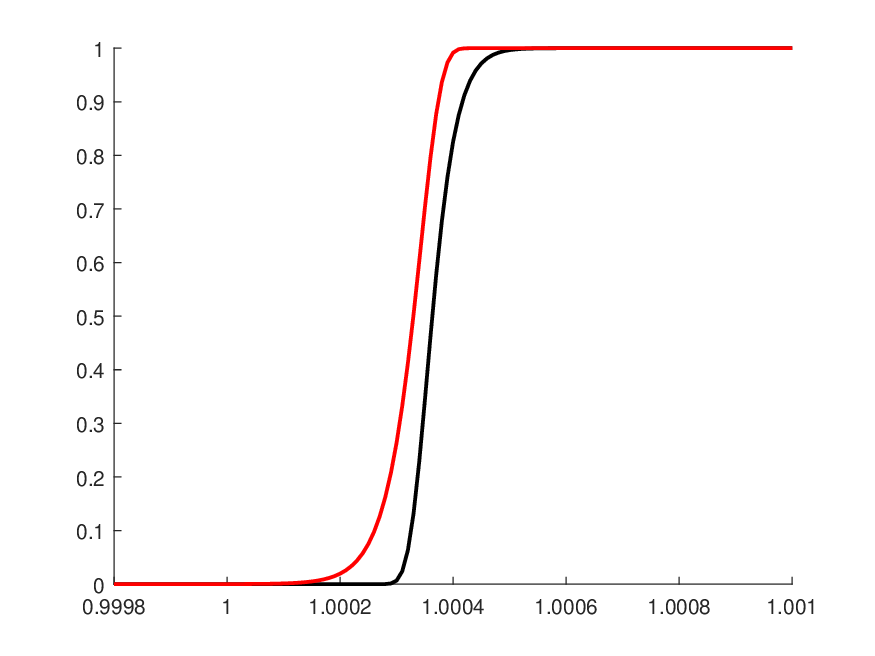}
\end{center}
	\caption{\footnotesize \textbf{(left panel)} The exact cumulative distribution function (cdf) of the scaled spectral radius $R_{N}$ for $N = 10^{3}$ (blue), $N = 10^{4}$ (blue dotted) and $N = 10^{7}$ (black) and $\alpha = 1$. The asymptotic cdf of $R_{N}$ (red dotted curve). \textbf{(right panel)} The exact formulation of the cumulative distribution function of the scaled spectral radius $R_{N}$ (black curve) and the asymptotic cumulative distribution function (red curve) for $N = 10^{7}$. The parameter $\alpha = 1$ for each curve.}
	\label{ComplexInducedGinibreScaledRmaxCdf}
\end{figure}

The scaled spectral radius and scaled minimum modulus of an $N \times N$ matrix from the complex induced Ginibre ensemble, with proportional rectangularity index, follow a Gumbel distribution in the limit as $N$ goes to infinity.
More precisely, the scaled spectral radius $R_{N} = \frac{r_{max}^{(N)}(G)}{\sqrt{L+N}}$ of an $N \times N$ matrix from the complex induced Ginibre ensemble with proportional rectangularity index, i.e., $L=\alpha N$ with $\alpha > 0$, follows a Gumbel distribution for maxima at the outer edge as $N$ goes to infinity. 
\\
\\
At the outer edge, setting $a = 1 + \frac{f_{N}(x)}{\sqrt{(1+\alpha)N}}$, where $f_{N}(x)$ is an increasing function in both $x$ and $N$, it is found that
\begin{equation}\label{ProbaSpectralRadiusRN}
	\mathcal{P}_N\left(a\right) 
	 	= P\left(\frac{r_{max}^{(N)}(G)}{\sqrt{(1+\alpha)N}} \leq a \right) 
	 	= \prod_{k=1}^{N} P\left( \frac{1}{(1+\alpha)N}\sum_{j=1}^{N+L-k+1}Z_{j} \leq a^{2} \right)
\end{equation}
where the $Z_{j}$ for $ j = \lbrace 1, \cdots, N+L-k+1 \rbrace$ are independent and identically distributed random variables following a standard exponential distribution. 
\\
\begin{theorem}\label{RMaxLimitTheoremComplexInducedGinibre} Let $G$ denote an $N \times N$ matrix from the complex induced Ginibre ensemble with rectangularity index $L$ proportional to $N$ (i.e., $L = \alpha N$, $\alpha > 0 $) and let $R_{N} = \frac{r_{max}^{(N)}(G)}{\sqrt{(1+\alpha)N}}$ denote the scaled spectral radius. \\ \\ Then,
\\ 
\begin{equation}
	 \lim_{N \rightarrow +\infty} P\left[ R_{N} \leq 1  + \sqrt{\frac{\gamma_{\alpha, N}}{2(1+\alpha)N}} +\frac{x}{2\sqrt{2(1+\alpha)N\gamma_{\alpha, N}}}\right]  
	 = e^{-e^{-x}}
\end{equation}\text{}\\ \text{}\\
where $\gamma_{\alpha, N} =  \log{\sqrt{(1+\alpha)N/2\pi}} - \log\log{N}$. \text{}\\ \text{}\\
The scaled spectral radius $R_{N}$ is approximated, as $N$ goes to infinity, by
\begin{equation}
	R_{N} \simeq  1 + T_{\alpha, N} + \xi_{\alpha, N}
\end{equation}
where $ T_{\alpha, N} = \sqrt{\frac{\gamma_{\alpha, N}}{2(1+\alpha)N}}$. The random variable $\xi_{\alpha, N} = \frac{X}{2\sqrt{2(1+\alpha)N\gamma_{\alpha, N}}}$, where $X$ follows a standard Gumbel distribution for maxima.\\
\end{theorem}\text{}\\
This result is similar to the one detailed in \cite{Rider2003} for the complex Ginibre ensemble. The exact formula of the probability density function $p_{R_{N}}\left( r \right)$ of the scaled spectral radius $R_{N}$, for a proportional rectangularity index $L = \alpha N $, $\alpha > 0$, is derived from equation (\ref{ProbaSpectralRadiusRN}) and is defined as follows for any value of $N$,
\\
\begin{align}
	p_{R_{N}}\left( r \right) 
		=  C^{(\alpha, N)}(r)\prod_{k=1}^{N}\frac{\gamma\left(k+\alpha N, (1+\alpha)Nr^{2}\right)}{\Gamma(k+\alpha N)} \sum_{j=1}^{N}\left[  \frac{ [(1+\alpha)N]^{j} r^{2j}}{\gamma \left(j+\alpha N, (1+\alpha)N r^{2}\right)}  \right]
\end{align}
where $\gamma(.,.)$ is the lower incomplete Gamma function and the function
\begin{equation}
	C^{(\alpha, N)}(r) = 2[(1+\alpha)N]^{\alpha N}r^{2(\alpha N-\frac{1}{2})}e^{-(1+\alpha)Nr^{2}} 
\end{equation}
\\
Numerical illustrations of these results are depicted in Figures \ref{ComplexInducedGinibreScaledSpectralRadiusRmaxProportionalL}, \ref{ComplexInducedGinibreLargeNRmaxDistribution} and \ref{ComplexInducedGinibreScaledRmaxCdf}. The analytical formulation of the probability density function $p_{R_{N}}$ of the scaled spectral radius $R_{N}$ fits exactly the empirical probability density function created from the generation of $10^{4}$ complex induced Ginibre matrices with proportional rectangularity index, as it is presented in Figure \ref{ComplexInducedGinibreScaledSpectralRadiusRmaxProportionalL}. As the size of these matrices increases, the probability density function of $R_{N}$ narrows close to the outer radius $r_{out} = 1$ (Figure \ref{ComplexInducedGinibreLargeNRmaxDistribution}). 

\text{}\\

\begin{figure}[H]
\begin{center}
	\centering		
	\includegraphics[scale=.4]{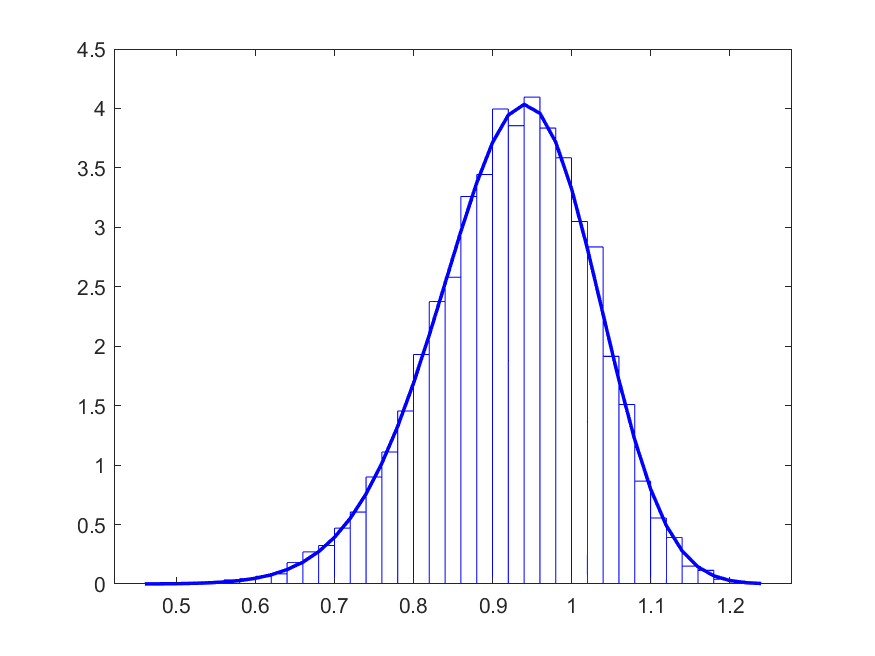}
	\includegraphics[scale=.4]{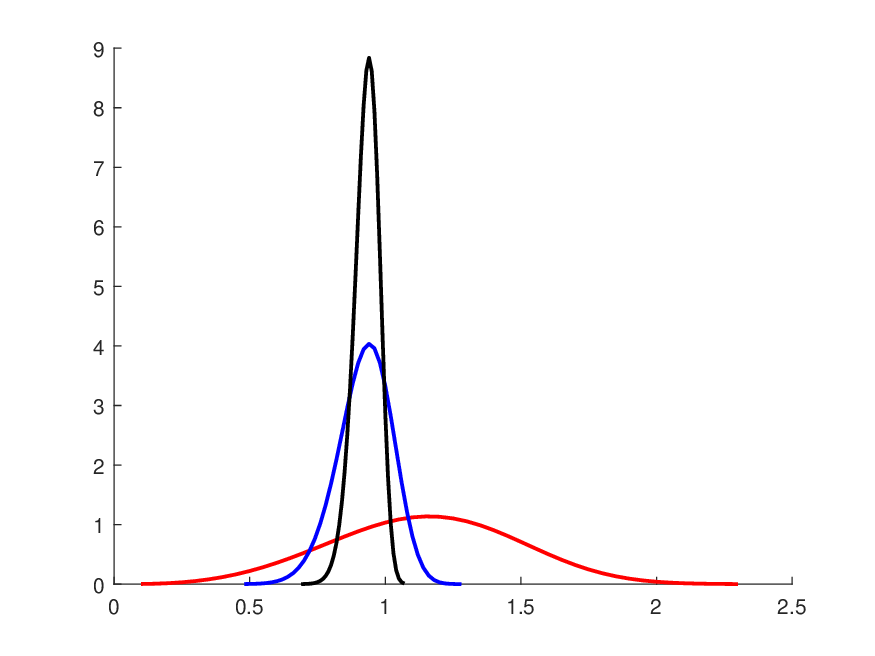}
\end{center}
\caption{\footnotesize \textbf{(left panel)} The empirical probability density function (histogram) of the scaled minimum modulus $r_{N}$ for $K=10^{4}$ generated $N \times N$ matrices from the complex induced Ginibre ensemble with $N =100$ and a proportional rectangularity index $L = \alpha N$, $\alpha = \frac{1}{10}$. The exact (analytical) probability density function is presented with the solid curve. \textbf{(right panel)} The analytical (exact) probability density function of the scaled minimum modulus $r_{N}$ for $K=10^{4}$ generated $N \times N$ matrices from the complex induced Ginibre ensemble with $N = 100$. The rectangularity index $L$ is proportional to $N$, i.e., $L = \alpha N$, $\alpha > 0$. The results are presented for different values of $\alpha = \left\lbrace \frac{1}{100}, \frac{1}{10}, \frac{4}{10} \right\rbrace$. }	\label{ComplexInducedGinibreScaledRminProportionalL}
\end{figure}

\begin{figure}[H]
	\begin{center}
		\includegraphics[scale=0.5]{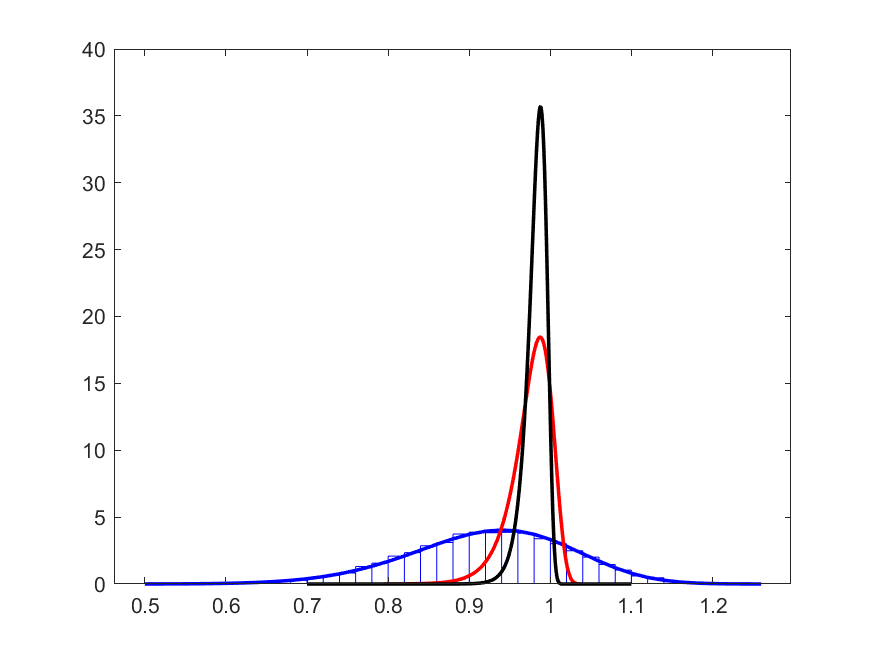}
	\end{center}	
	\caption{\footnotesize Empirical probability density function (histogram) and the exact (analytical) probability density function (blue curve) of the scaled minimum modulus $r_{N}$ for $K=10^{4}$ generated $N \times N$ matrices from the complex induced Ginibre ensemble with $N =100$ and a proportional rectangularity index $L = \alpha N$, $\alpha = \frac{1}{10}$. Limiting probability density functions of the scaled minimum modulus $r_{N}$ represented by the red and black curves for large $N = 10^{4}$ and $N = 2 \times 10^{4}$, respectively.}\label{ComplexInducedGinibreLargeNRminDistribution}
\end{figure}

\begin{figure}[H]
		\includegraphics[scale=.45]{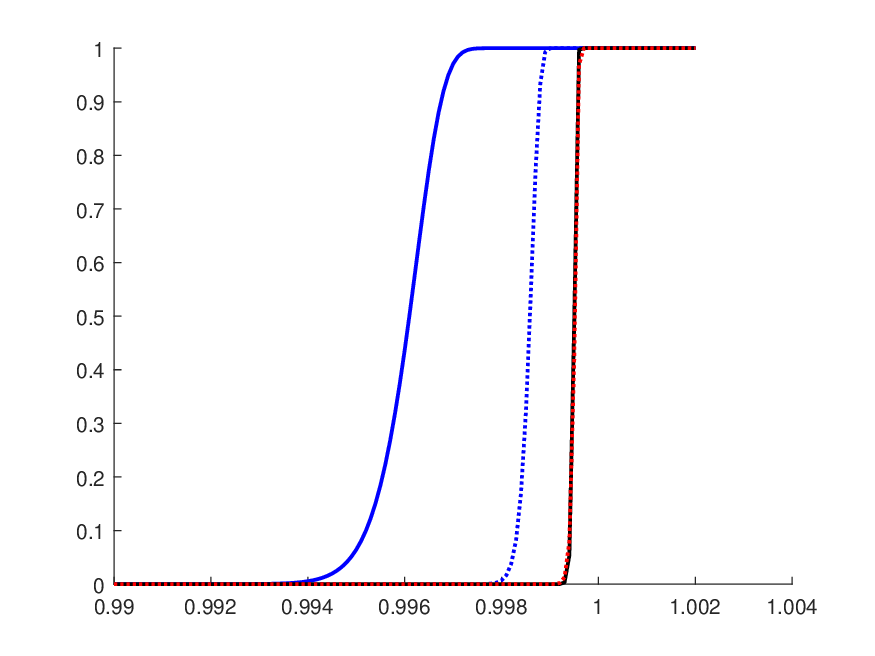}
	\hfill
		\includegraphics[scale=.45]{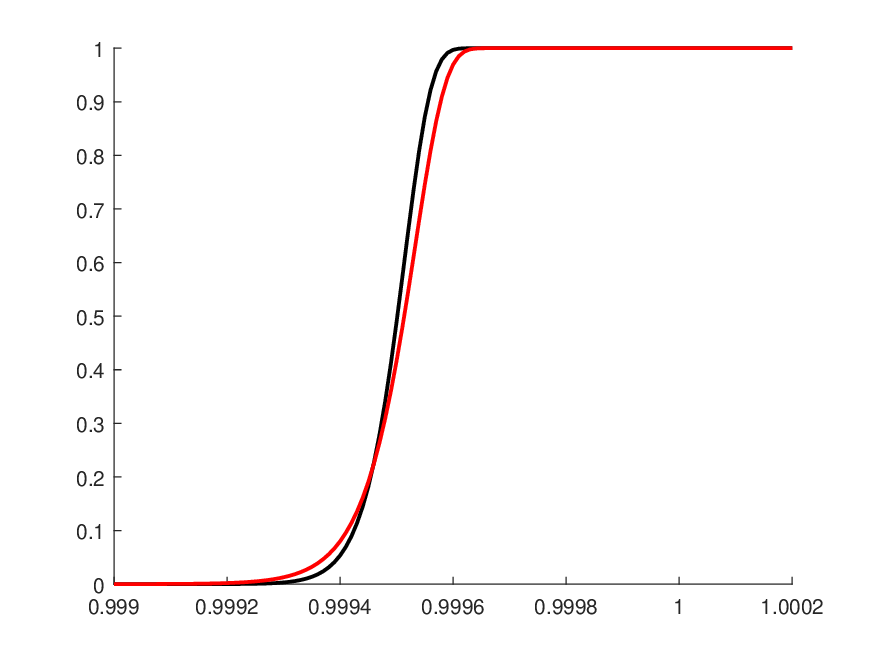}
\caption{\footnotesize \textbf{(left panel)} The exact cumulative distribution function (cdf) of the scaled minimum modulus $r_{N}$ for $N = 10^{5}$ (blue), $N = 10^{6}$ (blue dotted), $N = 10^{7}$ (black) and $\alpha = 1$. The asymptotic cdf of $r_{N}$ (red dotted curve). \textbf{(right panel)} The exact formulation of the scaled minimum modulus $r_{N}$ cumulative distribution function (black curve) and its asymptotic cumulative distribution function (red curve) for $N= 10^{7}$. The parameter $\alpha = 1$ for each curve.}
	\label{ComplexInducedGinibreScaledRminCdf}
\end{figure}

\text{}\\
The scaled minimum modulus of a matrix from the complex induced Ginibre ensemble with a proportional rectangularity index, follows a Gumbel distribution for minima at the inner edge (Theorem~\ref{RMinLimitTheoremComplexInducedGinibre}). Setting $a = 1 - \frac{f_{N}(x)}{\sqrt{\alpha N}}$ where $f_{N}(x)$ is an increasing function in both $x$ and $N$, it is found that
\begin{equation}\label{ProbaMinimumModulusrN}
	 \mathcal{P}_N\left( a \right) = P\left(\frac{r_{min}^{(N)}(A)}{\sqrt{\alpha N }} \geq a \right)
	   = \prod_{k=1}^{N} P\left( \frac{1}{\alpha N}\sum_{j=1}^{k+L}Z_{j} \geq a^{2} \right)
\end{equation}
where the $Z_{j}$ for $ j = \lbrace  1, \cdots, k+L  \rbrace$ are independent and identically distributed random variables, each following a standard exponential distribution.\\ 
\begin{theorem}\label{RMinLimitTheoremComplexInducedGinibre} Let $G$ denote an $N \times N$ matrix from the complex induced Ginibre ensemble with rectangularity index $L$ proportional to $N$ (i.e., $L = \alpha N$, $
\alpha > 0 $) and let $r_{N} = \frac{r_{min}^{(N)}(G)}{\sqrt{\alpha N}}$ denote the scaled minimum modulus.\\\\
Then, 
\begin{equation}
	 \lim_{N \rightarrow +\infty} P\left[ r_{N} \geq 1  - \sqrt{\frac{\gamma_{\alpha, N}}{2\alpha N}} +\frac{x}{2\sqrt{2 \alpha N\gamma_{\alpha, N}}}  \right]  = e^{-e^{x}}
\end{equation}\\ \\
where $\gamma_{\alpha, N} = \log{\sqrt{\alpha N/2\pi}} - \log\log{N}$. This limit corresponds to the survival distribution function of the standard Gumbel distribution for minima.\\ 
\\ 
The scaled minimum modulus $r_{N}$ is approximated as $N$ goes to infinity, by
\begin{equation}\label{AsymptoticScaledMinimumModulus}
	r_{N}  \simeq 1 - T_{\alpha, N} + \xi_{\alpha, N}
\end{equation}
where $ T_{\alpha, N} = \sqrt{\frac{\gamma_{\alpha, N}}{2 \alpha N}}$. The random variable $\xi_{\alpha, N} = \frac{X}{2\sqrt{2 \alpha N\gamma_{\alpha, N}}}$, where $X$ follows a standard Gumbel distribution for minima.
\end{theorem}\text{}\\
The exact formula of the probability density function $p_{r_{N}}\left(r \right)$ of the scaled minimum modulus $r_{N}$ for proportional rectangularity index $L = \alpha N$, $\alpha > 0,$ is derived from equation (\ref{ProbaMinimumModulusrN}) and is defined as follows for any value of $N$
\\
\begin{equation}
	p_{r_{N}}\left(r \right) = 2(\alpha N)^{\alpha N} r^{2(\alpha N-\frac{1}{2})} e^{-\alpha N r^{2}}\prod_{k=1}^{N}\frac{\Gamma(k+ \alpha N, \alpha Nr^{2})}{\Gamma(k+ \alpha N)}  \sum_{j=1}^{N}\left[  \frac{ (\alpha N)^{j} r^{2j} }{\Gamma(j+\alpha N, \alpha N r^{2})}   \right]
\end{equation}
where $\Gamma(.,.)$ is the upper incomplete Gamma function.\\
\\
Numerical results are presented in Figures \ref{ComplexInducedGinibreScaledRminProportionalL},  \ref{ComplexInducedGinibreLargeNRminDistribution} and \ref{ComplexInducedGinibreScaledRminCdf}. The analytical formulation of the probability density function of the scaled minimum modulus $r_{N}$ fits exactly its empirical counterpart from the generation of $10^{4}$ complex induced Ginibre matrices (left panel) with proportional rectangularity index $L = \alpha N$, $\alpha > 0$. As the size of these matrices increases, the probability density function of $r_{N}$ narrows around the inner radius $r_{in} = 1$. 

\text{}\\
The exact formula of the cumulative distribution function of $r_{N}$ converges towards the asymptotic distribution presented in Equation (\ref{AsymptoticScaledMinimumModulus}), as $N$ goes to infinity (Figure \ref{ComplexInducedGinibreScaledRminCdf}).\\
\text{}\\

\subsubsection{Independence of the spectral radius and the minimum modulus for the complex induced Ginibre ensemble}\label{IndependenceExtremeModuliComplexInducedGinibre}\text{}\\ \\
In this section, the independence of the spectral radius $r_{max}^{(N)}(G)$ and the minimum of moduli $r_{min}^{(N)}(G)$ is proved as $N$ goes to infinity, for an $N \times N$ matrix $G$ from the complex induced Ginibre ensemble.
\\
\begin{theorem}\label{ExtremeModuliIndependenceComplexGinibre} Let $\mathcal{B}_{L+N}$ and $\mathcal{B}_{L}$ denote Borel sets in the neighbourhood of $\sqrt{L+N}$ and $\sqrt{L}$, respectively. Let $r$ and $R$ denote reals in $\mathcal{B}_{L}$ and $\mathcal{B}_{L+N}$, respectively. The events $\lbrace r_{max}^{(N)}(G) \leq R \rbrace$ and $\lbrace r_{min}^{(N)}(G) \geq r \rbrace$ are independent in the limit as $N$ goes to infinity for fixed rectangularity index $L$. 
\end{theorem}
\begin{proof} Let $r_{max}^{(N)}(G)$ and $r_{min}^{(N)}(G)$ denote the spectral radius and the minimum modulus of an $N \times N$ complex induced Ginibre matrix $G$, respectively. The proof of Theorem~\ref{ExtremeModuliIndependenceComplexGinibre} is as follows. \\
\\
The probability
\begin{equation*}
	P\left( r_{min}^{(N)}(G) \geq r \text{  and  } r_{max}^{(N)}(G) \leq R \right) 
	=  \int_{r\leq \vert \lambda_{1}\vert \leq R} \cdots \int_{r\leq \vert \lambda_{N}\vert \leq R}  P_{N}(\lambda_{1}, \cdots, \lambda_{N}) \prod_{j=1}^{N}d^{2}\lambda_{j} 
\end{equation*}\\
where $r$ and $R$ are in $\mathcal{B}_{L}$ and $\mathcal{B}_{L+N}$, respectively. The joint probability density function $P_{N}(\lambda_{1}, \cdots, \lambda_{N})$ of the eigenvalues $ \lambda_{j}, j \in \lbrace 1, ..., N \rbrace$, for the complex induced Ginibre ensemble is presented in Equation~(\ref{ProbabilityDensityFunctionComplexInducedGinibre}). \\ 
\\
Let $r_{max}^{(N)}(G) \in \mathcal{B}_{L+N}$ and $r_{min}^{(N)}(G) \in \mathcal{B}_{L}$, where $\mathcal{B}_{L+N}$ is a Borel set in the neighbourhood of $\sqrt{L+N}$ and $\mathcal{B}_{L}$ a Borel set in the neighbourhood of $\sqrt{L}$. Applying Andreief's integration formula, the probability
\begin{equation*}
	P\left( r_{min}^{(N)}(G) \geq r \text{  and  } r_{max}^{(N)}(G) \leq R \right) 
	= P\left( r_{max}^{(N)}(G) \leq  R\right)  \prod_{k=1}^{N}\left[ 1 -  \frac{\gamma(k+L, r^{2})}{\gamma(k+L, R^{2})}\right]
\end{equation*}\text{}\\
For any $R $ in the Borel set $\mathcal{B}_{L+N}$, $\lim_{N \longrightarrow +\infty}\gamma(k+L, R^{2}) = \Gamma(k+L)$ for fixed rectangularity index $L > 0$. \\
\\
Thus,  
\begin{equation*}
	 \lim_{N \longrightarrow +\infty}\prod_{k=1}^{N}\left[ 1 -  \frac{\gamma(k+L, r^{2})}{\gamma(k+L, R^{2})}\right] 
=  \prod_{k=1}^{\infty} \frac{\Gamma(k+L, r^{2})}{\Gamma(k+L)} 
\end{equation*}\\
\\
and
\begin{equation*}
	 \lim_{N \rightarrow + \infty}P\left( r_{min}^{(N)}(G) \geq r \text{  and  } r_{max}^{(N)}(G) \leq R \right) 
	  = P\left( r_{max}^{(\infty)}(G) \leq R\right)P\left( r_{min}^{(\infty)}(G) \geq r\right)  
\end{equation*}
\end{proof}\text{}\\
\begin{theorem} Let $r_{max}^{(N)}(G)$ and $r_{min}^{(N)}(G)$ denote the spectral radius and the minimum modulus of an $N \times N$ matrix $G$ from the complex induced Ginibre ensemble, respectively. 
The scaled spectral radius $ \frac{r_{max}^{(N)}(G)}{\sqrt{L+N}}$ and the scaled minimum modulus  $ \frac{r_{min}^{(N)}(G)}{\sqrt{L+N}} $ are independent random variables, where for any $\alpha > 0$, $L = \alpha N$, in the limit as $N$ goes to infinity.\text{}\\
\end{theorem}

\begin{proof} Setting the rectangularity index as proportional to $N$, i.e., $L = \alpha N$, for any $\alpha > 0 $.
\begin{align*}
	& \lim_{N \rightarrow + \infty} P\left( \frac{r_{min}^{(N)}(G)}{\sqrt{L+N}} \geq r \text{  and  } \frac{r_{max}^{(N)}(G)}{\sqrt{L+N}} \leq R \right) 
	  =  \lim_{N \rightarrow + \infty}P\left( \frac{r_{max}^{(N)}(G)}{\sqrt{L+N}} \leq R\right) \lim_{N \rightarrow + \infty}  \prod_{k=1}^{N}\left[ 1 -  \frac{ \gamma(k+\alpha N, \left( 1 + \alpha \right)N r^{2})}{\gamma(k+\alpha N, \left( 1 + \alpha \right)NR^{2})}\right]
\end{align*}
\\
\begin{remark}\label{GammaFunctionRepresentation} 
\begin{equation*}
	\gamma(k + \alpha N, \left( 1 + \alpha \right)NR^{2}) = \Gamma(k + \alpha N) - \varepsilon_{k,\alpha N}(R)
\end{equation*}
where 
\begin{equation*}
	 \varepsilon_{k,N}^{(\alpha)}(R) = \Gamma(k + \alpha N, \left( 1 + \alpha \right)NR^{2})
\end{equation*}
and for any $k$ in $\mathbb{N}^{*}$,
\end{remark}\text{}\\
Then,
\begin{align*}
	& \lim_{N \rightarrow + \infty} \prod_{k=1}^{N}\left[ 1 - \frac{\gamma(k + \alpha N, \left( 1 + \alpha \right)Nr^{2})}{\gamma(k+ \alpha N, \left( 1 + \alpha \right)NR^{2})}\right]
	  = \lim_{N \rightarrow + \infty} \prod_{k=1}^{N}\left[ 1 - \frac{\gamma(k + \alpha N, \left( 1 + \alpha \right)Nr^{2})}{\Gamma(k + \alpha N) - \varepsilon_{k,N}^{(\alpha)}(R)}\right]\\
	& = \lim_{N \rightarrow + \infty} \prod_{k=1}^{N} \left(1 - \frac{\gamma(k + \alpha N, \left( 1 + \alpha \right)Nr^{2})}{\Gamma(k+ \alpha N)} \right)\lim_{N \rightarrow + \infty} \prod_{k=1}^{N} \left( 1 +  E_{k,N}^{(\alpha)}(R,r)   \right) 
\end{align*}
where $	E_{k,N}^{(\alpha)}(R,r) = \frac{\gamma(k+ \alpha N, \left( 1 + \alpha \right)Nr^{2})\varepsilon_{k,N}^{(\alpha)}(R)}{\Gamma(k+ \alpha N)^{2}}  + O\left( \frac{\gamma^{2}(k + \alpha N, \left( 1 + \alpha \right)Nr^{2})\varepsilon_{k,N}^{(\alpha)}(R)}{(\Gamma(k +\alpha N))^{3}}\right)$.\\ 
\\
\\ 
The limit of the partial product $ \prod_{k=1}^{N} \left( 1 +  E_{k,N}^{(\alpha)}(R,r) \right)$ is established from the limit of its lower and upper bounds, as $N$ goes to infinity.\\  \\
\\ 
\textit{\textbf{Lower bound.}} For any $r$ and $R\in \mathbb{R}^{+*}$, \\  
\begin{align*}
	\prod_{k=1}^{N} \left( 1 +  E_{k,N}^{(\alpha)}(R,r)  \right)	
	& > \prod_{k=1}^{N} \left( 1 + \frac{\gamma(k + \alpha N, \left( 1 + \alpha \right)Nr^{2})\varepsilon_{k,N}^{(\alpha)}(R)}{N\Gamma(k+ \alpha N)^{2}} \right) \\
	& > \prod_{k=1}^{N} \left( 1 + \frac{\gamma(1 + \alpha N, \left( 1 + \alpha \right)Nr^{2})\varepsilon_{k,N}^{(\alpha)}(R)}{N^{2}\Gamma^{2}(1+ \alpha N)} \right) 
	  = 1 + O\left( \frac{\gamma(1 + \alpha N, \left( 1 + \alpha \right)Nr^{2})\varepsilon_{k,N}^{(\alpha)}(R)}{N\Gamma^{2}(1+ \alpha N)} \right)
\end{align*}\\ 
\\
\textit{\textbf{Upper bound.}} Let $x=\left( 1 + \alpha \right)NR^{2}$. This implies, for any $R$ and $r \in \mathbb{R}^{+*}$,  
\begin{align*}
	\prod_{k=1}^{N} \left( 1 + E_{k,N}^{(\alpha)}(R,r)   \right)
	& = \prod_{k=1}^{N} \left( 1 +  \frac{\gamma(k+ \alpha N, \left( 1 + \alpha \right)Nr^{2})\varepsilon_{k,N}^{(\alpha)}(R)}{\Gamma(k+ \alpha N)^{2}}  + O\left( \frac{\gamma^{2}(k+ \alpha N, \left( 1 + \alpha \right)Nr^{2})\varepsilon_{k,N}^{(\alpha)}(R)}{(\Gamma(k+ \alpha N))^{3}}\right)\right)\\
	& < \prod_{k=1}^{N} \left( 1 +  O\left(\frac{\varepsilon_{k,N}^{(\alpha)}(R)}{\Gamma(k + \alpha N)} \right) \right)\\
	& = \exp \left[ \sum_{k=1}^{N}\log \left( 1  + O\left( \frac{\varepsilon_{k,N}^{(\alpha)}(R)}{\Gamma(k+\alpha N)}\right)\right) \right]\\
	& = \exp \left[ \sum_{k=1}^{N} O\left( \frac{\varepsilon_{k,N}^{(\alpha)}(R)}{\Gamma(k+ \alpha N)}\right) \right]
	  = \exp \left[ O\left(e^{-x^{2}}\sum_{j=0}^{k + \alpha N - 1}\frac{ x^{j+1}}{j!} \right) \right]
\end{align*} 
where $x$ goes to infinity, as $N$ goes to infinity. This implies that $e^{-x^{2}}\sum_{j=0}^{k+\alpha N - 1}\frac{ x^{j+1}}{j!}$ converge to zero.\\ \\
\\
\\
Applying the squeeze theorem, 
\begin{equation*}
	\lim_{N \rightarrow + \infty} \prod_{k=1}^{N} \left( 1 + E_{k,N}^{(\alpha)}(R,r)  \right) = 1
\end{equation*}
\\
Consequently, 
\begin{align*}
	\lim_{N \rightarrow + \infty} \prod_{k=1}^{N}\left[ 1 - \frac{\gamma(k+ \alpha N, \left( 1 + \alpha \right)r^{2})}{\gamma(k+ \alpha N, \left( 1 + \alpha \right)R^{2})}\right]
	& = \lim_{N \rightarrow + \infty} \prod_{k=1}^{N} \frac{\Gamma(k + \alpha N, \left( 1 + \alpha \right)r^{2})}{\Gamma(k + \alpha N)}
	  = \lim_{N \rightarrow + \infty} P \left( \frac{r_{min}^{(N)}(G)}{\sqrt{\left( 1 + \alpha \right)N}} \geq r  \right)
\end{align*} 
\end{proof} \text{}\\ 
\begin{theorem}[Independence of the scaled spectral radius and the scaled minimum modulus for the complex induced Ginibre ensemble]\label{IndependenceRMaxRminLimitTheoremComplexInducedGinibre} Let $G$ denote an $N \times N$ matrix from the complex induced Ginibre ensemble with proportional rectangularity index $L=\alpha N$, $\alpha = O\left( N^{m} \right)$, $m > 1$. The scaled spectral radius $R_{N}$ and the scaled minimum modulus $r_{N}$ are independent random variables, under scaling $\sqrt{\alpha N}$, as $N$ goes to infinity.\\ \\ 
More precisely, let $R_{N} = \frac{r_{max}^{(N)}(G)}{\sqrt{\alpha N}}$ and $r_{N} = \frac{r_{min}^{(N)}(G)}{\sqrt{\alpha N}}$. Setting $\rho = \sqrt{\frac{1+\alpha}{\alpha}}$, 
\begin{align*}
	 \lim_{N \rightarrow +\infty} P\left[ r_{N} \geq  1 - \sqrt{\frac{\gamma_{\alpha, N}}{2\alpha N}} + \xi_{\alpha}^{(N)}(y)
	 \text{}\text{ and }\text{} R_{N} \leq  \rho + \sqrt{\frac{\gamma_{\alpha, N}}{2\rho^{2}\alpha N}} + \eta_{\alpha}^{(N)}(x) \right]  
	  = e^{-e^{y}}e^{-e^{-x}}
\end{align*}
with $\gamma_{\alpha, N} = \log{\sqrt{\alpha N/2\pi}} - \log\log{N}$ and where
\begin{equation*}
	\xi_{\alpha}^{(N)}(y) = \frac{y}{2\sqrt{2 \alpha N\gamma_{\alpha, N}}} \text{  } \text{   and   } \text{  } \eta_{\alpha}^{(N)}(x) =  \frac{x}{2\rho\sqrt{2\alpha N\gamma_{\alpha, N}}}
\end{equation*}
The scaled spectral radius $R_{N}$ and the scaled minimum of moduli $r_{N}$ are approximated, as $N$ goes to infinity, by 
\begin{align*}
	R_{N} & \simeq \rho + \sqrt{\frac{\gamma_{\alpha, N}}{2\rho^{2}\alpha N}} + \frac{X}{2\rho\sqrt{2\alpha N\gamma_{\alpha, N}}}
\text{  }\text{  and  }\text{  }
	r_{N}  \simeq  1 - \sqrt{\frac{\gamma_{\alpha, N}}{2 \alpha N}} + \frac{Y}{2\sqrt{2 \alpha N\gamma_{\alpha, N}}} 
\end{align*}
\\
where $X$ and $Y$ follow a standard Gumbel (maximum) and standard Gumbel (minimum) distribution, respectively.
\end{theorem}\text{}\\
The proof of Theorem \ref{IndependenceRMaxRminLimitTheoremComplexInducedGinibre} is detailed in the Appendix \ref{IndependenceRMaxRminComplexInducedGinibre}.\\
\\
\subsection{Limiting left and right tails of the distribution of the minimum modulus for matrices from the complex Ginibre and complex induced Ginibre ensembles. A comparison.} \label{ComparisonsLimitingDistribution}\text{}\\ 
\\
The limiting left and right tails of the probability distribution of the minimum modulus for matrices from the complex Ginibre and complex induced Ginibre ensembles are investigated in this section. Derived in Section \ref{Methods}, the survival distribution function of the minimum modulus $r_{min}^{(N)}(A)$ of an $N \times N$ complex Ginibre matrix $A$, is used to determine the corresponding analytical probability density function
\begin{equation}
	p_{r_{min}^{(N)}(A)}\left(r \right) = 2r e^{-r^{2}} \prod_{k=0}^{N-1}\frac{\Gamma(k+1, r^{2})}{\Gamma(k+1)}\sum_{j=0}^{N-1}\left[  \frac{r^{2j}}{\Gamma(j+1, r^{2})} \right]
\end{equation}
\\
It represents the exact formulation of the probability density function of the minimum modulus $r_{min}^{(N)}(A)$.
\\
\begin{theorem}\label{LimitingLeftTailRminGinibre} Let $A$ denote an $N \times N$ random matrix from the complex Ginibre ensemble. The limiting left tail of the probability density function of the minimum modulus $r_{min}^{(N)}(A)$ is the probability density function of the Rayleigh distribution with parameter $\frac{1}{\sqrt{2}}$, as $N$ goes to infinity. More precisely, for $0 < r \ll 1$,
\begin{equation}
	\lim_{N \longrightarrow +\infty} p_{r_{min}^{(N)}(A)}\left(r\right) =  2re^{-r^{2}}\left( 1 + O(r^{2})  \right)
\end{equation}\\
\end{theorem}
\begin{proof} The limit of the partial product $\prod_{k=0}^{N-1}\frac{\Gamma(k+1, r^{2})}{\Gamma(k+1)}$, as $N$ goes to infinity, is studied in Section \ref{LimitTheoremComplexGinibre} for the complex Ginibre ensemble.\\
\\
For $0 < r \ll 1$,
\begin{equation*}
	\lim_{N \rightarrow +\infty}\prod_{k=0}^{N-1}\frac{\Gamma(k+1, r^{2})}{\Gamma(k+1)} 
	= e^{-r^{2}}\left( 1 - O\left(r^{4} \right) \right)
\end{equation*}
\\
The upper incomplete Gamma function has asymptotic expansion $\Gamma(k,x) \sim \Gamma(k) - \sum_{n=0}^{\infty}(-1)^{n}\frac{x^{k+n}}{n!(k+n)}$, as $x \longrightarrow 0^{+}$. \\
\\
The limit of the partial sum is derived as follows. Knowing that, for $0 < r \ll 1$,
\begin{equation*}
	\Gamma(j+1, r^{2}) = j!\left( 1 + O\left( r^{2} \right) \right)
\end{equation*} 
This implies,
\begin{align*}
	  \lim_{N \rightarrow +\infty}\sum_{j=0}^{N-1}\left[ \frac{1}{\Gamma(j+1, r^{2})}2r^{2j+1}e^{-r^{2}} \right]
	& = 2re^{-r^{2}}\sum_{j=0}^{+\infty}\frac{r^{2j}}{j!\left( 1 + O\left( r^{2} \right) \right)} 
	  = 2r \left( 1 + O(r^{2})\right)
\end{align*}\\
\\
Thus, the limiting left tail of the probability density function of the minimum modulus for an $N \times N$ matrix from the complex Ginibre ensemble is 
\begin{align*}
	\lim_{N \longrightarrow +\infty} p_{r_{min}^{(N)}(A)}\left(r\right) 
	& =  2re^{-r^{2}}\left( 1 + O(r^{2})\right)
\end{align*}\\
This function corresponds to the probability density function of the Rayleigh distribution with parameter $\sigma = \frac{1}{\sqrt{2}}$. \\
\end{proof}\text{}\\
This result is in line with the result acknowledged in Theorem \ref{RminComplexGinibreRayleighTheorem}.
\\ 
\\
Now, considering any rectangularity index $L \geq 0$, the limiting left tail of the distribution of the minimum modulus $r_{min}^{(N)}(G)$ for an $N \times N$ random matrix $G$ from the complex induced Ginibre ensemble, as $N$ goes to infinity, is derived as follows.
\\
\begin{theorem}\label{RminWeibull}
Let $G$ denote an $N \times N$ matrix from the complex induced Ginibre ensemble with rectangularity index $L \geq 0$. The left tail of the distribution of minimum modulus $r_{min}^{(N)}(G)$ is the Weibull distribution with shape parameter $\kappa = 2(L+1)$ and scale parameter $\lambda = \left( (L+1)! \right)^{1/\kappa}$ in the limit as $N$ goes to infinity. More precisely, as $r$ approaches zero,
\begin{equation}
	\lim_{N \longrightarrow +\infty}P\left( r_{min}^{(N)}(G) < r \right) =  1 - \exp\left[ -\frac{r^{2(L+1)}}{(L+1)!}  \right] + O\left( r^{2(L+2)} \right)
\end{equation}
\end{theorem}
\begin{proof} This result is derived using the same methodological approach as presented in section \ref{LimitingDistributionRmin}. \\ \\ For $0 < r \ll 1$, 
\begin{align*}
	\lim_{N \rightarrow +\infty}P\left( r_{min}^{(N)}(G) \geq r \right)  
	&= \sum_{j=0}^{+\infty} \frac{(-1)^{j}}{j!}\left[  \sum_{k=1}^{+\infty} \left(\frac{r^{2(k+L)}}{(k+L)!} + O\left( r^{2(k+L+1)} \right) \right) \right]^{j}
\end{align*}
\\
Consequently, in the limit as r goes to zero
\begin{equation*}
     \lim_{N \rightarrow +\infty} P\left( r_{min}^{(N)}(G) < r\right) 
	 	=  1 - \exp\left[ -\frac{r^{2(L+1)}}{(L+1)!} \right] + O\left( r^{2(L+2)} \right)
\end{equation*}
This corresponds to the Weibull distribution with shape parameter $\kappa = 2(L+1)$ and scale parameter \\ $\lambda = \left( (L+1)! \right)^{1/\kappa}$.
\end{proof}\text{}\\

\begin{corollary} Setting the rectangularity index $L$ to zero, the limiting left tail of the distribution of the minimum modulus $ r_{min}^{(N)}(G)$ corresponds to the limiting left tail of the distribution of the minimum modulus $ r_{min}^{(N)}(A)$, i.e., for $0< r \ll 1$, 
\begin{equation}
	\lim_{N \rightarrow +\infty} P\left( r_{min}^{(N)}(G) < r\right) =  1 - e^{ -r^{2} } + O(r^{4}) 
\end{equation}
\end{corollary}\text{}\\
\begin{lemma}\label{ProbabilityDensityFunctionRminComplexInducedGinibre}
The probability density function of the minimum modulus $r_{min}^{(N)}(G)$ of matrices from the complex induced Ginibre ensemble, with rectangularity index $L \geq 0$, is 
\begin{equation}
	p_{r_{min}^{(N)}(G)}\left(r\right) 
	=  2\frac{e^{-r^{2}}}{r}\prod_{k=1}^{N}\frac{\Gamma(k+L, r^{2})}{\Gamma(k+L)}\sum_{j=1}^{N}\left[  \frac{ r^{2(j+L)} }{\Gamma(j+L, r^{2})}   \right]
\end{equation}
\end{lemma}\text{}\\
\begin{corollary} For $0 < r \ll 1$, as $N$ goes to infinity, and for a rectangularity index $L=0$,
\begin{equation}
	\lim_{N \longrightarrow +\infty}p_{r_{min}^{(N)}(G)}\left( r\right)
		= \frac{\kappa}{\lambda}\left( \frac{r}{\lambda} \right)^{\kappa-1} e^{-(r/\lambda)^{\kappa}} \left( 1 + O(r^{2}) \right) 
\end{equation}
with Weibull distribution parameters $\kappa=2$ and $\lambda = 1$. This corresponds to the probability density function of the Rayleigh distribution with parameter $\sigma = \frac{1}{\sqrt{2}}$.\\
\end{corollary}
\begin{proof}Let $G$ denote a $N \times N$ complex induced Ginibre matrix with rectangularity index $L=0$. \\ 
\\
For $0 < r \ll 1$, from Lemma~\ref{ProbabilityDensityFunctionRminComplexInducedGinibre} and the methodological approach used to prove Theorem~\ref{LimitingLeftTailRminGinibre},
\begin{align*}
	\lim_{N \longrightarrow +\infty}p_{r_{min}^{(N)}(G)}\left( r\right)
	& =  \frac{2e^{-r^{2}}}{r}\prod_{k=0}^{N-1}\frac{\Gamma(k+1, r^{2})}{\Gamma(k+1)}\sum_{j=0}^{N-1}\left[ \frac{r^{2(j+1)}}{\Gamma(j+1, r^{2})} \right] \\
	& = \frac{\kappa}{\lambda}\left( \frac{r}{\lambda} \right)^{\kappa-1} e^{-(r/\lambda)^{\kappa}} \left( 1 + O(r^{2}) \right)\\
	& = \lim_{N \longrightarrow +\infty} p_{r_{min}^{(N)}(A)}\left( r\right) 
\end{align*}
with Weibull distribution parameters $\kappa=2$ and $\lambda = 1$.\\ 
\end{proof}
\text{}\\
The limiting right tail of the distribution of the minimum modulus for the complex induced Ginibre ensemble, in the limit as $N$ goes to infinity for fixed rectangularity index $L \geq 0$, is derived using the Euler-Maclaurin summation formula and the asymptotic expansion of the upper incomplete Gamma function as $r$ goes to infinity.\\
\\
Thus,
\begin{equation*}
	 \lim_{N \rightarrow +\infty} P\left( r_{min}^{(N)}(G) \geq r\right) 
	 	= \exp{\left[ \sum_{k=1}^{+\infty}\log\left( \frac{\Gamma(k+L, r^{2})}{\Gamma(k+L)} \right)	\right]}
\end{equation*}\\ \\
Now, the asymptotic expansion of the upper incomplete Gamma function (cf. \cite{Temme1995}) is
\begin{equation*}
	\Gamma(a, z) = z^{a-1}e^{-z} + O(z^{a-2}e^{-z}) \text{, as } z \text{ goes to infinity}
\end{equation*} 
Thus, for large $r$,
\begin{equation*}
	\frac{1}{2}\log\left( \frac{\Gamma(L+1, r^{2})}{\Gamma(L+1)} \right) 
	= \frac{1}{2}\left[ -r^{2} + \log\left(\frac{r^{2L} \left(1 + O\left(  \frac{1}{r^{2}}\right)\right)}{\Gamma(L+1)}  \right)\right]
\end{equation*}\\
\\ 
A precise evaluation of the integral $\int_{1}^{+\infty}\log\left(\frac{\Gamma(k+L, r^{2})}{\Gamma(k+L)}\right)dk$ is $ -\frac{r^{4}}{4} + O(r^{2})$.\\ \\
\\
With the use of the Euler-Maclaurin summation formula, the infinite series 
\begin{equation*}
	\sum_{k=1}^{+\infty}\log\left( \frac{\Gamma(k+L, r^{2})}{\Gamma(k+L)} \right) = \int_{1}^{+\infty}\log\left( \frac{\Gamma(k+L, r^{2})}{\Gamma(k+L)}\right)dk + \frac{1}{2}\log\left( \frac{\Gamma(L+1, r^{2})}{\Gamma(L+1)}\right) + O\left( r \right)
\end{equation*}\\
\\
Consequently,
\begin{equation*}
	\sum_{k=1}^{+\infty}\log\left( \frac{\Gamma(k+L, r^{2})}{\Gamma(k+L)} \right) 
		= -\frac{r^{4}}{4} + \frac{1}{2}\left[ -r^{2} + \log\left(\frac{r^{2L} \left(1 + O\left(  \frac{1}{r^{2}}\right)\right)}{\Gamma(L+1)}  \right)\right] + O\left( r^{2} \right)
\end{equation*}\\
\\
This implies that, for large $r$, 
\begin{equation}
	\lim_{N \rightarrow +\infty} P\left( r_{min}^{(N)}(G) < r\right) 
		= 1 - e^{-\frac{r^{4}}{4}}\exp\left( -\frac{r^{2}}{2} + \frac{1}{2}\log\left(\frac{  r^{2L} + r^{2(L-1)}O\left( 1  \right) }{ \Gamma(L+1)} \right) + O\left( r^{2} \right)\right)
\end{equation}\\
\\
For a rectangularity index $L=0$, for large $r$, and an $N \times N$ complex Ginibre matrix $A$, \\
\begin{equation}
	 \lim_{N \rightarrow +\infty} P\left( r_{min}^{(N)}(G) < r\right)  
	 	= 1 - e^{-\frac{r^{4}}{4}\left( 1 + O\left( \frac{1}{r^{2}} \right) \right)} 
	 	= \lim_{N \rightarrow +\infty}P\left( r_{min}^{(N)}(A) < r\right) 
\end{equation} \\ 
\\
\section{Conclusions}\label{Conclusions}\text{}\\
This research was conducted from pure academic interest in the fields of Random Matrix Theory and Probability theory. The survival distribution function of the minimum modulus of an $N \times N$ matrix from the complex Ginibre ensemble is analytically derived with the use of Andreief's integration formula \cite{Forrester2018}. It corresponds to an $N$-th partial product of regularised upper incomplete Gamma functions as in references \cite{GrobeHaakeSommers1988, Forrester1992, AkemannPhillipsShifrin2009}.\\ 
The scaled extreme moduli of an $N \times N$ complex induced Ginibre matrix follow a Gumbel distribution as N goes to infinity. More precisely, their asymptotic distribution is well approximated as a shifted and scaled standard Gumbel distribution.
The stochastic modelling of these random variables might reveal more accurate results while taking into account their mean and second central moment within a similar methodological framework as the one established in reference \cite{Rider2003}. The analytical formulations of the latter statistical features might have relevance.\\ 
Derived for the non-Hermitian ensembles considered in this paper, the analytical formulations of the probability density functions of the extreme moduli (the spectral radius and the minimum modulus with and without scaling) exactly coincide with their empirical counterpart sampled from thousands of matrices. 
The independence of the minimum modulus and the spectral radius is a trivially accepted fact for which no formal proof was detailed for the complex Ginibre ensemble, as well as for its random matrix ensemble generalisation. The independence of these random variables has been proved in the present paper for these two random matrix ensembles at appropriate scalings. It is acknowledged to hold as $N$ goes to infinity.\\ 
The limiting left and right tails of the distribution of the minimum modulus have also been investigated for large sizes of complex Ginibre matrices. They are the Rayleigh distribution with parameter $\sigma = 1 / \sqrt{2}$ and the Weibull distribution with shape parameter $\kappa = 4$ and scale parameter $\lambda = \kappa^{1/\kappa}$, respectively. \\
Limiting left and right tails of the distribution of the minimum modulus provide insight into the exact nature of the distribution of this extreme modulus for the complex Ginibre ensemble, in the limit as $N$ goes to infinity. The limiting left and right tails of distribution of the minimum modulus of a matrix from the complex Ginibre ensemble are the same as the limiting left and right tails of distribution of the minimum modulus of a matrix from the complex induced Ginibre ensemble for a rectangularity index equal to zero. They are the Rayleigh and the Weibull distributions, respectively. \\
The limiting distribution of the minimum modulus for the complex Ginibre ensemble, as $N$ goes to infinity, is absent from the literature. A mixture of distributions involving either the Rayleigh distribution or the Weibull distribution might define the right statistics of this minimum. This is left for forthcoming research.\\
\newpage
\begin{appendices}
\section{Proofs of limit theorems for the complex induced Ginibre ensemble}\label{moduliExtremaLimitingDistributionComplexInducedGinibre}\text{}\\
\subsection{Proof of Theorem 6}\label{AppendixRmaxTheorem}\text{}\\
\text{}\\
Let $G$ denote an $N \times N $ complex induced Ginibre matrix. The limiting distribution of the scaled spectral radius is determined at the outer edge of the ring, in the limit as $N$ goes to infinity and with a rectangularity index $L$ proportional to $N$, i.e., $L = \alpha N$, $\alpha > 0$.\text{}\\
\text{}\\
Let $r_{out} = \sqrt{(1+\alpha)N}$ denote the outer radius. The chosen scaling is $\sqrt{\left(1+\alpha\right) N}$. Applying Andreief's integration formula \cite{Forrester2018}, the cumulative distribution function of the scaled spectral radius $\frac{r_{max}^{(N)}(G)}{\sqrt{\left(1 + \alpha\right)N}}$, is \text{}\\
\text{}\\
\begin{align*}
	   P\left(\frac{r_{max}^{(N)}(G)}{\sqrt{(1+\alpha)N}} \leq a \right) 
	 & = \frac{1}{\pi^{N}\prod_{k=1}^{N}\Gamma(k + L)} \prod_{k=0}^{N-1}   \pi\int_{0}^{(1+\alpha)N a^{2}} t^{k+L} e^{-t} dt\\
	 & = \frac{1}{\prod_{k=1}^{N}\Gamma(k + L)} \prod_{k=0}^{N-1} \int_{0}^{(1+\alpha)N a^{2}} t^{k+L} e^{-t} dt\\
	 & = \prod_{k=1}^{N} \int_{0}^{(1+\alpha)N a^{2}} \frac{ t^{(k+L-1)} e^{-t}}{\Gamma(k + L)} dt\\
	 & = \prod_{k=1}^{N}\int_{0}^{(1+\alpha)Na^{2}}f_{Gamma(k+L,1)}(t)dt 
	   = \prod_{k=1}^{N} \frac{\gamma(k + L,(1+\alpha)N a^{2})}{\Gamma(k + L)}
\end{align*}\\
The function $f_{Gamma(k+L, 1)}(t)$ is the probability density function of the Gamma distribution with shape parameter $k+ L$ and rate parameter $1$.\text{}\\
\text{}\\
This implies that,
\begin{align*}
	P\left(\frac{r_{max}^{(N)}(G)}{\sqrt{(1+\alpha)N}} \leq a \right)
	& = \prod_{k=1}^{N} P\left( \frac{1}{(1+\alpha)N}\sum_{j=1}^{N+L-k+1}Z_{j} \leq a^{2} \right)
\end{align*} 
\text{}\\
where the $Z_{j}$ for $ j = \lbrace 1, \cdots, N-k+L+1 \rbrace$ are independent and identically distributed random variables, each following an exponential distribution with parameter $1$. \\ 
\\
\begin{align*}
	\mathcal{P}_{N}\left( 1 + \frac{f_{N}(x)}{\sqrt{(1+\alpha)N}} \right)	
	&  = \prod_{k=1}^{N} P\left( \frac{1}{(1+\alpha)N}\sum_{j=1}^{N+L-k+1}Z_{j} \leq a^{2} \right)\\
	& = \prod_{k=1}^{N} P\left( \frac{1}{\sqrt{(1+\alpha)N}}\sum_{j=1}^{N+L-k+1}(Z_{j}-1) \leq \phi_{N}(x) + \frac{k-1}{\sqrt{(1+\alpha)N}}  \right) 
	  = \prod_{k=0}^{N-1}p_{k}
\end{align*}
where
\begin{equation*}
	p_{k} = P\left( \frac{1}{\sqrt{(1+\alpha)N}}\sum_{j=1}^{(1+\alpha)N-k}(Z_{j}-1) \leq \phi_{N}(x) + \frac{k}{\sqrt{(1+\alpha)N}}  \right) 
\end{equation*}\\
\\
Let $a = \frac{r_{out}}{\sqrt{(1+\alpha)N}} + \frac{f_{N}(x)}{\sqrt{(1+\alpha)N}} = 1 + \frac{f_{N}(x)}{\sqrt{(1+\alpha)N}}$ and $\phi_{N}(x) = 2f_{N}(x) + \frac{f_{N}^{2}(x)}{\sqrt{(1+\alpha)N}} = 2f_{N}(x)\left( 1 + \frac{f_{N}(x)}{2\sqrt{(1+\alpha)N}} \right)$.  The function $f_{N}(x) = o(\sqrt{(1+\alpha)N})$, meaning that $\lim_{N \rightarrow +\infty} \frac{f_{N}(x)}{\sqrt{(1+\alpha)N}} = 0$. \text{}\\ \text{}\\ The probability $\mathcal{P}_{N}\left(1 + \frac{f_{N}(x)}{\sqrt{(1+\alpha)N}}\right)$ is bounded by any partial product of $p_{k}$, i.e., $\prod_{k=0}^{(1+ \alpha)N \delta_{N}} p_{k}$, for a positive $\delta_{N}$ less than 1 (cf. \cite{Rider2003}). \text{}\\ \text{}\\ More precisely, \\
\begin{equation*}
	 \prod_{k=0}^{(1+ \alpha)N \delta_{N}} p_{k} \prod_{k=(1+ \alpha)N \delta_{N}}^{N-1} p_{k}\leq \mathcal{P}_{N}\left( 1 + \frac{f_{N}(x)}{\sqrt{(1+\alpha)N}}\right) \leq \prod_{k=0}^{(1 + \alpha)N \delta_{N}} p_{k}
\end{equation*}\\
\\
The following is derived using the same arguments as stated in \cite{Rider2003} and applying the Markov inequality as well as the definition of the quantile of the Exponential distribution. The function $\phi_{N}(x)$ is a positive and increasing function. This implies that
\begin{align*}
	&\prod_{k = (1+\alpha)N \delta_{N}}^{N-1} P\left( \frac{1}{\sqrt{(1+\alpha)N}}\sum_{j=1}^{(1+\alpha)N-k}(Z_{j}-1) \leq \phi_{N}(x) + \frac{k}{\sqrt{(1+\alpha)N}}  \right) \\
	& = \prod_{k = (1+\alpha)N\delta_{N}}^{N-1} P\left( \sum_{j=1}^{(1+\alpha)N-k}(Z_{j}-1) \leq \sqrt{(1+\alpha)N}\phi_{N}(x) + k \right) 
\end{align*}\\
\\
The random variables $Z_{j}$, $i \in \lbrace 1, \cdots, (1+\alpha)N-k \rbrace$ are independent and identically distributed. Applying the Markov inequality, with $0 < \eta < 1$, and as the exponential is an strictly increasing and convex function, this implies that\\
\\
\begin{align*}
	&\prod_{k = (1+\alpha)N \delta_{N}}^{N-1} P\left( \frac{1}{\sqrt{(1+\alpha)N}}\sum_{j=1}^{(1+\alpha)N-k}(Z_{j}-1) \leq \phi_{N}(x) + \frac{k}{\sqrt{(1+\alpha)N}}  \right) \\
	& = \prod_{k = (1+\alpha)N\delta_{N}}^{N-1} \left [ 1 - P\left( \sum_{j=1}^{(1+\alpha)N-k}(Z_{j}-1) > \sqrt{(1+\alpha)N}\phi_{N}(x) + k \right)\right] \\
	& =  \prod_{k = (1+\alpha)N\delta_{N}}^{N-1} P\left( \sum_{j=1}^{(1+\alpha)N-k}(Z_{j}-1) \leq \sqrt{(1+\alpha)N}\phi_{N}(x) + k \right) \\
	& =  \prod_{k = (1+\alpha)N\delta_{N}}^{N-1} P\left( \sum_{j=1}^{(1+\alpha)N-k}Z_{j} \leq \sqrt{(1+\alpha)N}\phi_{N}(x) + (1+\alpha)N \right) \\
	& \geq  \prod_{k = (1+\alpha)N\delta_{N}}^{N-1} P\left( \sum_{j=1}^{(1+\alpha)N-k}Z_{j} \leq \sqrt{(1+\alpha)N}\phi_{N}(x) \right)  \\
	& = \prod_{k = (1+\alpha)N\delta_{N}}^{N-1}  \left[ 1 -  P\left( \sum_{j=1}^{(1+\alpha)N-k}Z_{j} > \sqrt{(1+\alpha)N}\phi_{N}(x) \right)\right]\\
	& \geq  \prod_{k = (1+\alpha)N\delta_{N}}^{N-1} \left[  1 - e^{-\eta\sqrt{(1+\alpha)N}\phi_{N}(x)}E\left[ e^{ \eta Z_{1}} \right]^{(1+\alpha)N-k} \right]
\end{align*}\text{}\\ 
\text{}\\
Furthermore, \\
\begin{align*}
	e^{-\eta (1+\alpha)N (\frac{k}{(1+\alpha)N} -1) - (1+\alpha)N\ln(1-\eta)} 
	& \geq e^{((1+\alpha)N -k)(\eta - \ln (1 - \eta))}\\
	& = e^{-\eta (1+\alpha)N (\frac{k}{(1+\alpha)N} -1) - \left((1+\alpha)N-k \right)\ln(1-\eta)} > 1\\ 
\end{align*}
And with $0 < \eta < 1$ and $\forall \text{ } Y \in \mathbb{R}$,
\begin{equation*}
	e^{Y - \eta \sqrt{(1+\alpha)N}\phi_{N}(x)} < e^{Y} \Rightarrow 1-e^{Y - \eta \sqrt{(1+\alpha)N}\phi_{N}(x)}  > 1- e^{Y}
\end{equation*}\\
Now, setting $\eta = 1-\frac{1}{\frac{k}{(1+\alpha)N} - 1}$ which does maximise the remaining exponent (similarly presented in \cite{Rider2003}) and applying the quantile formula of the Exponential distribution, this implies\\
\\
\begin{align*}
	&\prod_{k = (1+\alpha)N \delta_{N}}^{N-1} P\left( \frac{1}{\sqrt{(1+\alpha)N}}\sum_{j=1}^{(1+\alpha)N-k}(Z_{j}-1) \leq \phi_{N}(x) + \frac{k}{\sqrt{(1+\alpha)N}}  \right) \\
	& \geq \prod_{k = (1+\alpha)N \delta_{N}}^{N-1} \left( 1 - e^{-(1+\alpha)N \left[ \eta \left(\frac{k}{(1+\alpha)N} -1 \right) + \ln(1-\eta )\right]} \right)\\
	& \geq \prod_{k = (1+\alpha)N \delta_{N}}^{N-1} \left( 1 - e^{-(1+\alpha)N \left[  \frac{k}{(1+\alpha)N} - \ln\left( \frac{k}{(1+\alpha)N} -1 \right)\right]} \right)\\
	& \geq \left( 1- e^{-(1+\alpha)N\delta_{N}^{2}} \right)^{N}
\end{align*}
\\
Finally,
\begin{equation*}
	\prod_{k=(1+ \alpha)N \delta_{N}}^{N-1} p_{k} \geq \left(1 - e^{-(1+\alpha)N\delta_{N}^{2})} \right)^{N}
\end{equation*} \\
The parameter $\delta_{N}$ is chosen such that $\left(1 - e^{-(1+\alpha)N\delta_{N}^{2})} \right)^{N} = \left( 1 - \frac{1}{N^{2}} \right)^{N}= 1 - O\left( \frac{1}{N} \right)$, i.e., $\delta_{N} = \sqrt{\frac{2\log{N}}{(1+\alpha)N}}$. \text{}\\ 
\text{}\\ 
Consequently,
\begin{equation*}
	\lim_{N \rightarrow +\infty}  \left[ \prod_{k=0}^{(1+ \alpha)N \delta_{N}} p_{k} \prod_{k=(1+ \alpha)N \delta_{N}}^{N-1} p_{k} \right] = \lim_{N \rightarrow +\infty} \prod_{k=0}^{(1+ \alpha)N \delta_{N}} p_{k} 
\end{equation*} \text{}\\
\text{}\\
Now, applying the squeeze theorem, uniformly in $N$ and $x$, for bounded $x > - \infty$,  
\begin{align*}
	\lim_{N \rightarrow +\infty}\log \left(\mathcal{P}_{N}\left(1 + \frac{f_{N}(x)}{\sqrt{(1+\alpha)N}} \right) \right) 
	& = \lim_{N \rightarrow +\infty} \log \left( \prod_{k=0}^{ \sqrt{2(1+\alpha)N\log(N)}} p_{k} \right) 
\end{align*}
Also,
\begin{align*}	
	& \prod_{k=0}^{(1+ \alpha)N \delta_{N}} P\left( \frac{1}{\sqrt{(1+\alpha)N}}\sum_{j=1}^{(1+\alpha)N-k}(Z_{j}-1) \leq \phi_{N}(x) + \frac{k}{\sqrt{(1+\alpha)N}}  \right) \\
	& = \prod_{k=0}^{(1+ \alpha)N \delta_{N}} P\left( \frac{1}{\sqrt{N+L-k}}\sum_{j=1}^{N+L-k}(Z_{j}-1) \leq \left(\sqrt{\frac{(1+\alpha)N}{N+L-k}} \right)\left( \phi_{N}(x) + \frac{k}{\sqrt{(1+\alpha)N}} \right)  \right) \\
	& = \prod_{k=0}^{(1+ \alpha)N \delta_{N}} p_{k}
\end{align*}
where 
\begin{equation*}
	p_{k} = P\left( \frac{1}{\sqrt{N+L-k}}\sum_{j=1}^{N+L-k}(Z_{j}-1) \leq \left( \sqrt{\frac{(1+\alpha)N}{N+L-k}} \right)\left( \phi_{N}(x) + \frac{k}{\sqrt{(1+\alpha)N}} \right)  \right)
\end{equation*}
\\
and, for fixed $k \in \left[ 0, (1+ \alpha)N \delta_{N} \right]$, the factor $\sqrt{\frac{(1+\alpha)N}{N+L-k}}$ goes to one as $N$ goes to infinity.\\ 
\\
Finally,
\begin{align*}
	& \lim_{N \rightarrow +\infty}\log \left(\mathcal{P}_{N}\left(1 + \frac{f_{N}(x)}{\sqrt{(1+\alpha)N}}\right) \right) \\
	& = \lim_{N \rightarrow +\infty} \sum_{k=0}^{\sqrt{2(1+\alpha)N \log{N}}}   \log P\left( \frac{1}{\sqrt{(1+\alpha)N-k}}\sum_{j=1}^{(1+\alpha)N-k}(Z_{j}-1) \leq \phi_{N}(x) + \frac{k}{\sqrt{(1+\alpha)N}}  \right)
\end{align*}
where $\delta_{N} = \sqrt{\frac{2\log{N}}{(1+\alpha)N}}$.\text{}\\
\text{}\\
The Edgeworth expansion is then used to get the probability density of the following standardised random variable  
\begin{align*}
	\frac{\sqrt{(1+\alpha)N-k}(\bar{Z}- \mu)}{\sigma} 
	& = \frac{1}{\sqrt{(1+\alpha)N-k}}\sum_{j=1}^{(1+\alpha)N-k}(Z_{j}-1)
\end{align*} 
where each element of sequence $Z_{1}, \cdots, Z_{(1+\alpha)N-k}$ is i.i.d. exponentially distributed with parameter equal to $1$. The random variable $\bar{Z} = \frac{1}{(1+\alpha)N-k}\sum_{j=1}^{(1+\alpha)N-k}Z_{j}$ is the empirical mean. The mean of $\bar{Z}$ is equal to $1$ and its standard deviation is $\frac{\sigma}{\sqrt{(1+\alpha)N-k}} = \frac{1}{\sqrt{(1+\alpha)N-k}}$ with $\sigma = 1$.\\ \\
\\
Applying the Edgeworth expansion, the logarithm of the probability $p_{k}$ is then
\begin{align*}
	& \log \left(  P\left( \frac{1}{\sqrt{(1+\alpha)N-k}}\sum_{j=1}^{(1+\alpha)N-k}(Z_{j}-1) \leq \phi_{N}(x) + \frac{k}{\sqrt{(1+\alpha)N}}  \right) \right) \\
	& = \log 
	\begin{pmatrix}
	\int_{-K_{N}}^{\phi_{N}(x) + \frac{k}{\sqrt{(1+\alpha)N}} } \frac{e^{-\frac{t^{2}}{2}}}{\sqrt{2\pi}}dt + O\left(\frac{1}{\sqrt{(1+\alpha)N-k}}\sup_{\vert c \vert  \leq  Y}  \phi_{N}^{2}(c)e^{-\frac{\phi_{N}^{2}(x)}{2}}\right)\\
	 + O\left(\frac{1}{(1+\alpha)N-k}\right) + O\left( K_{N}((1+\alpha)N-k)^{-3/2}\right) + O\left(e^{-\frac{K_{N}^{2}}{2}}\right)
	\end{pmatrix}\\
	& = \log \left(	
	\int_{-K_{N}}^{\phi_{N}(x) + \frac{k}{\sqrt{(1+\alpha)N}} } \frac{e^{-\frac{t^{2}}{2}}}{\sqrt{2\pi}}dt  \right) + O\left(\frac{1}{\sqrt{(1+\alpha)N-k}}\sup_{\vert c \vert  \leq  Y}  \phi_{N}^{2}(c)e^{-\frac{\phi_{N}^{2}(c)}{2}}\right)\\
	& + O\left(\frac{1}{(1+\alpha)N-k}\right) + O\left( K_{N}((1+\alpha)N-k)^{-3/2}\right) + O\left(e^{-\frac{K_{N}^{2}}{2}}\right)
\end{align*}
where, as in \cite{Rider2003}, for $x$ restricted as in $ \vert x \vert \leq Y$ for some large positive $Y$, $K_{N}$ goes to $+\infty$ faster than $\sup_{\vert x \vert \leq Y}\phi_{N}(x)$.\\
\begin{remark} Let
\begin{equation*}
	T_{N} = \int_{-K_{N}}^{\phi_{N}(x) + \frac{k}{\sqrt{(1+\alpha)N}} } \frac{e^{-\frac{t^{2}}{2}}}{\sqrt{2\pi}}dt 
\end{equation*}
For $z$ in the neighbourhood of zero, it is known that the Taylor expansion of  $\log(1+z) = z + O(z^{2})$
and $\lim_{N \longrightarrow +\infty} T_{N} = 1$.\text{}\\
\text{}\\
Then, in the limit of large values of $N$,
\begin{align*}
	&\log 
	\begin{pmatrix}
	T_{N} + O\left(\frac{1}{\sqrt{(1+\alpha)N-k}}\sup_{\vert c \vert  \leq  Y}  \phi_{N}^{2}(c)e^{-\frac{\phi_{N}^{2}(x)}{2}}\right)	 + O\left(\frac{1}{(1+\alpha)N-k}\right)\\ + O\left( K_{N}((1+\alpha)N-k)^{-3/2}\right) + O\left(e^{-\frac{K_{N}^{2}}{2}}\right)
	\end{pmatrix}
	\\
	& = \log 
	\left[ T_{N} 
	\begin{pmatrix}
	1 + \frac{1}{T_{N}}	
	\begin{pmatrix}
	 O\left(\frac{1}{\sqrt{(1+\alpha)N-k}}\sup_{\vert c \vert  \leq  Y}  \phi_{N}^{2}(c)e^{-\frac{\phi_{N}^{2}(x)}{2}}\right)\\
	 + O\left(\frac{1}{(1+\alpha)N-k}\right) + O\left( K_{N}((1+\alpha)N-k)^{-3/2}\right) + O\left(e^{-\frac{K_{N}^{2}}{2}}\right) 
	 \end{pmatrix}
	 \end{pmatrix}
	\right]\\
	& = \log  
	\left[ 
	1 + 
	\begin{pmatrix}
	 O\left(\frac{1}{\sqrt{(1+\alpha)N-k}}\sup_{\vert c \vert  \leq  Y}  \phi_{N}^{2}(c)e^{-\frac{\phi_{N}^{2}(x)}{2}}\right)\\
	 + O\left(\frac{1}{(1+\alpha)N-k}\right) + O\left( K_{N}((1+\alpha)N-k)^{-3/2}\right) + O\left(e^{-\frac{K_{N}^{2}}{2}}\right) 
	 \end{pmatrix} 	\right]\\
	& = \log\left[ \int_{-K_{N}}^{\phi_{N}(x) + \frac{k}{\sqrt{(1+\alpha)N}} } \frac{e^{-\frac{t^{2}}{2}}}{\sqrt{2\pi}}dt \right] \\
	& + \log 
	\left[
	1 +
	\begin{pmatrix}
	 O\left(\frac{1}{\sqrt{(1+\alpha)N-k}}\sup_{\vert c \vert  \leq  Y}  \phi_{N}^{2}(c)e^{-\frac{\phi_{N}^{2}(x)}{2}}\right)\\
	 + O\left(\frac{1}{(1+\alpha)N-k}\right) + O\left( K_{N}((1+\alpha)N-k)^{-3/2}\right) + O\left(e^{-\frac{K_{N}^{2}}{2}}\right) 
	 \end{pmatrix}
	\right]\\
	& = \log \left[	
	\int_{-K_{N}}^{\phi_{N}(x) + \frac{k}{\sqrt{(1+\alpha)N}} } \frac{e^{-\frac{t^{2}}{2}}}{\sqrt{2\pi}}dt  \right] + O\left(\frac{1}{\sqrt{(1+\alpha)N}}\sup_{\vert c \vert  \leq  Y}  \phi_{N}^{2}(c)e^{-\frac{\phi_{N}^{2}(x)}{2}}\right)\\
	& + O\left(\frac{1}{(1+\alpha)N}\right) + O\left( K_{N}((1+\alpha)N)^{-3/2}\right) + O\left(e^{-\frac{K_{N}^{2}}{2}}\right)
\end{align*}
\end{remark} \text{}\\
As mentioned in \cite{Rider2003}, the lower limit of integration in the leading term of the logarithm of the probability $p_{k}$ is extended from $-K_{N}$ down to $-\infty$. Defining the function $f(k)$ as
\begin{equation*}
	f(k) =  \log \left[ \int_{-\infty}^{\phi_{N}(x) + \frac{k}{\sqrt{(1+\alpha)N}}} \frac{e^{- \frac{t^{2}}{2}}}{\sqrt{2\pi}}dt \right]  
\end{equation*}
\begin{align*}
	\sum_{k=0}^{\sqrt{2(1+\alpha)N\log{N}}}f(k) 
	& = \int_{0}^{\sqrt{2(1+\alpha)N\log{N}}} \log \left[ \int_{-\infty}^{\phi_{N}(x) + \frac{t}{\sqrt{(1+\alpha)N}}} \frac{e^{- \frac{s^{2}}{2}}}{\sqrt{2\pi}}ds \right]dt 
	  + E_{N}
\end{align*}
where $E_{N}$ is an error term.\text{}\\ \text{}\\
Using the change of variables $u =  \phi_{N}(x) + \frac{t}{\sqrt{(1+\alpha)N}}$, this implies that $du = \frac{dt}{\sqrt{(1+\alpha)N}}$  and  $dt = \sqrt{(1+\alpha)N}du$
\\
\begin{equation*}
	\sum_{k=0}^{\sqrt{2(1+\alpha)N\log N}} f(k) = \sqrt{(1+\alpha)N} \int_{\phi_{N}(x) }^{ \phi_{N}(x)+ \sqrt{2\log{N}}}\log \left[ \int_{-\infty}^{u} \frac{e^{- \frac{s^{2}}{2}}}{\sqrt{2\pi}}ds\right]du + E_{N}
\end{equation*}
\\
where, with respect to the error terms  $O\left(\frac{1}{\sqrt{(1+\alpha)N-k}}\sup_{\vert c \vert  \leq  Y}  \phi_{N}^{2}(c)e^{-\frac{\phi_{N}^{2}(x)}{2}}\right)$ and $O\left(\frac{1}{(1+\alpha)N}\right)$ of the Edgeworth expansion of the logarithm of the probability $p_{k}$,
\begin{equation*}
	E_{N} = O\left(  \left(\frac{\sqrt{\log N}}{\sqrt{(1+\alpha)N}} \right) \vee  \left( \sqrt{\log{N}}\sup_{\vert c \vert \leq Y}\phi_{N}^{2}(c)e^{-\frac{\phi_{N}^{2}(c)}{2}} \right) \right)
\end{equation*}
\\
\begin{remark} Identification of the error term corresponding to the sum of the error of integration at each integration point $\phi_{N}(x) + \frac{k}{\sqrt{(1+\alpha)N}}$. \text{}\\
\text{}\\
With $\phi_{N}(x)= o\left(\sqrt{\log \left((1+\alpha)N \right)}\right)$ and $K_{N} = O(\log{(1+\alpha)N})$, this implies that\\
\begin{align*}
	&\left| \sqrt{(1+\alpha)N} \sum_{k=0}^{\sqrt{2(1+\alpha)N\log N}} \int_{\phi_{N}(x) + \frac{k}{\sqrt{(1+\alpha)N}}}^{\phi_{N}(x) + \frac{k+1}{\sqrt{(1+\alpha)N}}} \log \left( 1 + \int_{\phi_{N}(x) + \frac{k}{\sqrt{(1+\alpha)N}}}^{s} e^{-\frac{t^{2}}{2}}dt\right)ds  \right| \\
	& \leq C \left| \sqrt{\log{N}}e^{-\frac{\phi_{N}^{2}(x)}{2}} \right| 
		=  C \left| \sqrt{\log{N}}\left((1+\alpha)N\right)^{-1/2} \right| 
\end{align*}\\
where $C$ is a positive constant.\text{}\\
\text{}\\
This implies that
\\
\begin{equation*}
	\sqrt{(1+\alpha)N} \sum_{k=0}^{\sqrt{2(1+\alpha)N\log N}} \int_{\phi_{N} + \frac{k}{\sqrt{(1+\alpha)N}}}^{\phi_{N} + \frac{k+1}{\sqrt{(1+\alpha)N}}}\log \left(1 + 	\int_{\phi_{N} + \frac{k}{\sqrt{(1+\alpha)N}}}^{s} e^{-\frac{t^{2}}{2}}dt\right)ds 
	 = O\left( \frac{\sqrt{\log(N)}}{\sqrt{(1+\alpha)N}}  \right)
\end{equation*}\\ \\
This completes the remark.
\end{remark}\text{}\\ \\
The term $\frac{1}{\sqrt{(1+\alpha)N}}$ appearing in the lower limit of integration is negligible. Extending the upper limit to $+\infty$, the function $\phi_{N}(x) = o \left( \sqrt{\log \left( (1+\alpha)N  \right)}\right)$.\\ 
\\
This implies the following result
\begin{equation*}
	\log{\mathcal{P}_{N}}\left( 1 + \frac{f_{N}(x)}{\sqrt{(1+\alpha)N}} \right) = \sqrt{(1+\alpha)N}\int_{\phi_{N}(x)}^{+\infty}\log \left[ \int_{-\infty}^{t} \frac{e^{-\frac{s^{2}}{2}}}{2\pi}ds\right]dt + E_{N}
\end{equation*}
\\
Let $F_{\infty}(x)$ denote the limiting cumulative distribution function of the scaled spectral radius of an $N \times N$ complex induced Ginibre matrix $G$ as $N$ goes to infinity. The limit of its logarithm is
\begin{equation*}
	 \log{F_{\infty}}(x) = \lim_{N \rightarrow +\infty} \sqrt{(1+\alpha)N} \int_{\phi_{N}(x)}^{+\infty} \log \left[ 1 - \int_{t}^{+\infty} e^{-\frac{s^{2}}{2}}\frac{ds}{\sqrt{2\pi}} \right]dt
\end{equation*} 
where $\lim_{N \rightarrow +\infty} E_{N} = 0$.\\ \\
As the function $f_{N}$ goes to infinity, the logarithm of the limiting distribution $F_{\infty}(x)$ is 
\begin{align*}
	\log{F_{\infty}}(x) 
		& = -\lim_{N \rightarrow +\infty} \sqrt{\frac{(1+\alpha)N}{2\pi}} \int_{\phi_{N}(x)}^{+\infty}\left[
	 	\int_{t}^{+\infty}e^{-\frac{s^{2}}{2}}ds + O\left( \left( \int_{t}^{+\infty} e^{-\frac{s^{2}}{2}}ds \right)^{2}\right) \right]dt \\
	 	& =  -\lim_{N \rightarrow +\infty} \sqrt{\frac{(1+\alpha)N}{2\pi}} \int_{\phi_{N}(x)}^{+\infty}\left[
	 	\frac{1}{t}e^{-t^{2}/2} + O\left( \frac{e^{-t^{2}}}{t^{2}} \right) 	\right]dt 	 	
\end{align*}
\\
Furthermore, $\phi_{N}(x) = 2f_{N}(x) + \frac{f_{N}^{2}(x)}{\sqrt{(1+\alpha)N}} = 2f_{N}(x)\left( 1 + \frac{f_{N}(x)}{2\sqrt{(1+\alpha)N}} \right)$. The function $f_{N}(x) = o(\sqrt{(1+\alpha)N})$, meaning that $\lim_{N \rightarrow +\infty} \frac{f_{N}(x)}{\sqrt{(1+\alpha)N}} = 0$.\text{}\\
\text{}\\
This implies that 
\begin{align*}
	\lim_{N \rightarrow +\infty} \phi_{N}^{2}(x) 
	& = \lim_{N \rightarrow +\infty} 4f_{N}^{2}(x)\left( 1 + \frac{f_{N}(x)}{2\sqrt{(1+\alpha)N}} \right)^{2}\\
	& = \lim_{N \rightarrow +\infty} 4f_{N}^{2}(x)\left( 1 + \frac{f_{N}(x)}{\sqrt{(1+\alpha)N}} + \frac{f_{N}^{2}(x)}{4(1+\alpha)N} \right)
	  = \lim_{N \rightarrow +\infty} 4f_{N}^{2}(x)
\end{align*}
\\
Also,
\begin{equation*}
	\int_{\phi_{N}(x)}^{+\infty} \frac{1}{t}e^{-t^{2}/2} dt  =  \frac{1}{\phi_{N}^{2}(x)}e^{-\frac{\phi_{N}^{2}(x)}{2}} - 2\int_{\phi_{N}(x)}^{+\infty}\frac{e^{-t^{2}/2}}{t^{3}}dt
\end{equation*}
\\
and in the limit as $N$ goes to infinity,
\begin{equation*}
	\int_{\phi_{N}(x)}^{+\infty}\frac{e^{-t^{2}/2}}{t^{3}}dt = \frac{1}{\phi_{N}^{2}(x)}e^{-\frac{\phi_{N}^{2}(x)}{2}} \times \frac{1}{\phi_N^{2}(x)} + O\left( \int_{\phi_{N}(x)}^{+\infty} 	
	\frac{e^{-t^{2}/2}}{t^{5}}dt \right)
\end{equation*}
This implies that
\begin{align*}
	\log{F_{\infty}}(x) 		 	
	 & = -\lim_{N \rightarrow +\infty} \sqrt{\frac{(1+\alpha)N}{2\pi}} \left[\frac{1}{\phi_{N}^{2}(x)}e^{-\frac{\phi_{N}^{2}(x)}{2}} \left(1 + O\left(\frac{1}{\phi_{N}^{2}(x)}\right)\right) \right] \\
	 & = -\lim_{N \rightarrow +\infty} \sqrt{\frac{(1+\alpha)N}{2\pi}} \left[
	 	\frac{1}{4f_{N}^{2}(x)}e^{-2f_{N}^{2}(x)} \left(1 + O\left(\frac{1}{4f_{N}^{2}(x)}\right)\right) \right]
\end{align*}
\\
Then, a convenient choice of the function $f_{N}(x)$ is performed, such that the limiting cumulative distribution function $F_{\infty}(x)$ is from the class of Extreme Value Distributions composed of the three types of extreme value distributions for maxima.\\
\\
More precisely, choosing 
\begin{equation*}
	f_{N}^{2}(x) = \frac{1}{2}\log \left( \frac{e^{x}\sqrt{(1+\alpha)N/2\pi}}{\log N} \right)
\end{equation*}
\\
The limit, uniformly on compact sets in $x$, as in \cite{Rider2003},
\begin{align*}
	& \lim_{N \rightarrow +\infty} \log \mathcal{P}_{N}\left[ R_{N} \leq 1 + \sqrt{\frac{1}{2(1+\alpha)N}}\left( \log \frac{\sqrt{(1+\alpha)N/2\pi}}{\log{N}} + x \right)^{1/2}\right] \\
	& = -\lim_{N \rightarrow +\infty} \sqrt{\frac{(1+\alpha)N}{2\pi}} \left[
	 	\frac{1}{4f_{N}^{2}(x)}e^{-2f_{N}^{2}(x)} \left(1 + O\left(\frac{1}{4f_{N}^{2}(x)}\right)\right) \right]\\
	& = -\lim_{N \rightarrow +\infty}\exp(-x)\frac{\log{N}}{2\log(e^{x}\sqrt{(1+\alpha)N/(2\pi)} \times \frac{1}{\log{N}})}
\end{align*}
\\	
Also, 
\begin{equation*}
	  \frac{\log{N}}{2\log(e^{x}\sqrt{(1+\alpha)N/(2\pi)} \times \frac{1}{\log{N}})} \\
	   = \frac{1}{2}\frac{\log{N}}{\left[ x + \frac{1}{2}\log((1+\alpha)N/(2\pi)) - \log\log N) \right]}
\end{equation*}\\
\\
And, for any $\alpha > 0$ and fixed $\vert x \vert \ll N$,
\begin{align*}
	 & \lim_{N \rightarrow + \infty} \frac{1}{2}\frac{\log{N}}{\left[ x + \frac{1}{2}\log((1+\alpha)N/(2\pi)) - \log\log N) \right]} \\
	 & = \lim_{N \rightarrow + \infty}\frac{1}{2}\frac{\log{N}}{\left[\frac{1}{2}\log((1+\alpha)N/(2\pi)) - \log\log N \right]}\\
	 & = \lim_{N \rightarrow + \infty} \frac{\log{N}}{\log((1+\alpha)N/(2\pi))}\\
	 & = \lim_{N \rightarrow + \infty}\frac{1}{1 + \frac{\log  \left(\frac{1+\alpha}{2\pi} \right)}{\log N}}
\end{align*}
\\
This implies
\begin{equation*}
	\lim_{N \rightarrow + \infty} \frac{\log{N}}{2\log(e^{x}\sqrt{(1+\alpha)N/2\pi} 
	\times \frac{1}{\log{N}})} = 1
\end{equation*}\\
\\
Finally,  
\begin{equation*}
	\lim_{N \rightarrow +\infty} \log \mathcal{P}_{N}\left[ R_{N} \leq 1 + \sqrt{\frac{1}{2(1+\alpha)N}}\left( \log \frac{\sqrt{(1+\alpha)N/2\pi}}{\log{N}} + x \right)^{1/2}\right] = -\exp(-x)
\end{equation*}\\
\\
This limit is the logarithm of the cumulative distribution function of the standard Gumbel distribution for maxima.\text{}\\
\text{}\\
Furthermore, setting $\gamma_{\alpha, N} =  \log \frac{\sqrt{(1+\alpha)N/2\pi}}{\log{N}} = \log{\sqrt{(1+\alpha)N/2\pi}} - \log\log{N}$ and using the Taylor expansion of the square root function\\
\begin{equation*}
	\left( \log \frac{\sqrt{(1+\alpha)N/2\pi}}{\log{N}} + x \right)^{1/2} 
		=  \gamma_{\alpha, N}^{1/2} \left( 1 + \frac{x}{\gamma_{\alpha, N}} \right)^{1/2}
		= \gamma_{\alpha, N}^{1/2} + \frac{x}{2\gamma_{\alpha, N}^{1/2}} + O\left(\left(\frac{x}{\gamma_{\alpha, N}}\right)^{2}\right)
\end{equation*}\\
\\
This implies,
\begin{equation*}
	\lim_{N \rightarrow +\infty} \mathcal{P}_{N}\left[ R_{N} \leq 1  + \sqrt{\frac{\gamma_{\alpha, N}}{2(1+\alpha)N}} +\frac{x}{2\sqrt{2(1+\alpha)N\gamma_{\alpha, N}}}\right] = \exp(-\exp(-x))
\end{equation*}\\
\\
This limit corresponds to cumulative distribution function of the standard Gumbel distribution for maxima.\\ 
\\
The scaled spectral radius $R_{N}$ is a random variable approximated by
\begin{align*}
	R_{N} & \simeq  1 + \sqrt{\frac{\gamma_{\alpha, N}}{2(1+\alpha)N}} +\frac{X}{2\sqrt{2(1+\alpha)N\gamma_{\alpha, N}}} \\
		    & = 1 + \sqrt{\frac{\gamma_{\alpha, N}}{2(1+\alpha)N}} -\frac{\log(Z)}{2\sqrt{2(1+\alpha)N\gamma_{\alpha, N}}} \\
		    & = 1 + T_{\alpha, N} + \xi_{\alpha, N}
\end{align*}
where the term $ T_{\alpha, N} = \sqrt{\frac{\gamma_{\alpha, N}}{2(1+\alpha)N}}$ with $\gamma_{\alpha N} = \log{\sqrt{(1+\alpha)N/2\pi}} - \log\log{N}$. The random variable $Z$ denote a random variable following a standard Exponential distribution. The random variable $X = -\log(Z)$ has a standard Gumbel distribution for maxima. The random variable $\xi_{\alpha, N} = -\frac{\log(Z)}{2\sqrt{2(1+\alpha)N\gamma_{\alpha, N}}}$ has a Gumbel distribution. \\ 
\\  
This completes the proof of Theorem \ref{RMaxLimitTheoremComplexInducedGinibre}. \\
\\

\subsection{Proof of Theorem 7}\label{AppendixRminTheorem}\text{}\\
\text{}\\ 
Let $G$ denote an $N \times N $ complex induced Ginibre matrix. The limiting distribution of the scaled minimum modulus at the inner edge of the ring, in the limit as $N$ goes to infinity and with a rectangularity index $L$ proportional to $N$, i.e., $L = \alpha N$, $\alpha > 0$, is studied in this section.\text{}\\
\text{}\\
Let $r_{in} = \sqrt{\alpha N} $ denote the inner radius. The chosen scaling is $\sqrt{\alpha N}$. The survival distribution function of the scaled minimum modulus $\frac{r_{min}^{(N)}(G)}{\sqrt{\alpha N }}$ is \\
\begin{align*}
	 P\left(\frac{r_{min}^{(N)}(G)}{\sqrt{\alpha N }} \geq a \right) 
	 & = \frac{1}{\prod_{k=1}^{N}\Gamma(k + \alpha N)} \prod_{k=0}^{N-1} \int_{\alpha Na^{2}}^{+\infty} t^{k+L} e^{-t} dt\\
	 & = \prod_{k=1}^{N}\int_{\alpha Na^{2}}^{+\infty}\frac{ t^{k+L -1} e^{-t}}{\Gamma(k + L)}dt
	   = \prod_{k=1}^{N}\int_{\alpha Na^{2}}^{+\infty}f_{Gamma(k+L, 1)}(t)dt 
	   = \prod_{k=1}^{N}\frac{\Gamma\left(k + \alpha N , \alpha Na^{2} \right)}{\Gamma(k + \alpha N)}
\end{align*}\\
\\
The function $f_{Gamma(k+L, 1)}(t)$ is the probability density function of the Gamma distribution with shape parameter $k+ L$ and rate parameter $1$.\\ 
\\
This implies,
\begin{equation*}
	P\left(\frac{r_{min}^{(N)}(G)}{\sqrt{\alpha N}} \geq a \right)
	 = \prod_{k=1}^{N} P\left( \sum_{j=1}^{k+L}Z_{j} \geq \alpha N a^{2} \right)
   	 = \prod_{k=1}^{N} P\left( \frac{1}{\alpha N}\sum_{j=1}^{k+L}Z_{j} \geq a^{2} \right)
\end{equation*} 
where the $Z_{j}$ for $ j = \lbrace  1, \cdots, k+L  \rbrace$ are independent and identically distributed random variables following an Exponential distribution with parameter $1$.\text{}\\
\text{}\\
Let $a = 1 - \frac{f_{N}(x)}{\sqrt{\alpha N}}$, with $\alpha > 0$. The function $f_{N}(x) = o\left( \sqrt{\alpha N} \right)$, meaning that $\lim_{N \rightarrow +\infty}\frac{f_{N}(x)}{\sqrt{\alpha N}} = 0$. \\ 
\\
Then,
\begin{align*}
\mathcal{P}_{N}\left( 1 - \frac{f_{N}(x)}{\sqrt{\alpha N}} \right) 
	& = \mathcal{P}_{N}\left(\frac{\min_{1 \leq k \leq N} \vert z_{k} \vert}{\sqrt{\alpha N}} \geq 1 -  \frac{f_{N}(x)}{\sqrt{\alpha N}} \right)
	  = \prod_{k=1}^{N} p_{k}
\end{align*}
where 
\begin{equation*}
	p_{k} =  P\left( \frac{1}{\sqrt{\alpha N}}\sum_{j=1}^{k+L}\left( Z_{j} - 1 \right) \geq -\phi_{N}(x) - \frac{k}{\sqrt{\alpha N}} \right)
\end{equation*}\\
\\
with $\phi_{N}(x) = 2f_{N}(x)\left(1- \frac{f_{N}(x)}{2\sqrt{\alpha N}}\right)$. \\ \\
\\
The probability $\mathcal{P}_{N}(1 - \frac{f_{N}(x)}{\sqrt{\alpha N}})$ is bounded by any partial product of $p_{k}$, i.e., bounded by the partial product $\prod_{k=0}^{\alpha N \delta_{N}} p_{k}$, for a positive $\delta_{N}$ less than 1 (cf. \cite{Rider2003}). Let $\delta_{N} = \sqrt{\frac{2\log N}{\alpha N}}$.\\
\\
Then,  
\begin{equation*}
	 \prod_{k=1}^{k=\alpha N \delta_{N}} p_{k}\prod_{\alpha N \delta_{N}}^{N} p_{k}  \leq \mathcal{P}_{N}\left( 1 - \frac{f_{N}(x)}{\sqrt{\alpha N}}\right) \leq \prod_{k=1}^{\alpha N \delta_{N}} p_{k}
\end{equation*}\\
\\
The random variable $ \frac{1}{\sqrt{\alpha N+k}}\sum_{j=1}^{k+L}\left( Z_{j} -1 \right)$ converges in distribution to a standard normal variable as $N$ goes to infinity (The Central Limit Theorem).\text{}\\
\text{}\\
Thus, 
\begin{align*}
	& \prod_{k=\alpha N \delta_{N}}^{N}  P\left( \frac{1}{\sqrt{\alpha N}}\sum_{j=1}^{k+L}\left( Z_{j} -1 \right) \geq -\phi_{N}(x) - \frac{k}{\sqrt{\alpha N}} \right)\\
	& = \prod_{k=\alpha N \delta_{N}}^{N}  P\left( \frac{1}{\sqrt{\alpha N + k }}\sum_{j=1}^{\alpha N + k}\left( Z_{j} -1 \right) \geq -\frac{\sqrt{\alpha N}}{\sqrt{\alpha N +k}}\phi_{N}(x) - \frac{k}{\sqrt{\alpha N +k}} \right)\\
	& > \prod_{k=\alpha N \delta_{N}}^{N}  P\left( \frac{1}{\sqrt{\alpha N + k }}\sum_{j=1}^{\alpha N + k}\left( Z_{j} -1 \right) \geq - \frac{k}{\sqrt{\alpha N +k}} \right)\\
	& =  \prod_{k=\alpha N \delta_{N}}^{N}  \left[ 1 - \frac{1}{2}\left( 1 + \erf\left(-\frac{\left( \frac{k}{\sqrt{\alpha N +k}}\right)}{\sqrt{2}}\right) \right) \right]\\
	& > \prod_{k=\alpha N \delta_{N}}^{N} \left( 1 - e^{-\frac{k^{2}}{\alpha N}} \right)
	> \prod_{k=\alpha N \delta_{N}}^{N} \left( 1 - e^{-\alpha N \delta_{N}^{2}} \right)
	> \left( 1 - \frac{1}{N^{2}} \right)^{N}= 1 - O\left(\frac{1}{N}\right)\\
\end{align*}
Consequently,
\begin{equation*}
	\lim_{N \rightarrow +\infty}  \left[ \prod_{k=1}^{\alpha N \delta_{N}} p_{k} \prod_{k=\alpha N \delta_{N}}^{N} p_{k} \right] = \lim_{N \rightarrow +\infty} \prod_{k=1}^{\alpha N \delta_{N}} p_{k} 
\end{equation*}\text{}\\
\text{}\\
Now, applying the squeeze theorem, uniformly in $N$ and $x$, for bounded $x > - \infty$,  
\begin{align*}
	&\lim_{N \rightarrow +\infty}\log \left(\mathcal{P}_{N}\left(1 - \frac{f_{N}(x)}{\sqrt{\alpha N}}\right) \right) 
	  = \lim_{N \rightarrow +\infty} \log \left( \prod_{k=1}^{\alpha N \delta_{N}} p_{k} \right) 
	  = \lim_{N \rightarrow +\infty} \log \left( \prod_{k=1}^{\sqrt{2\alpha N\log(N)}} p_{k} \right)
\end{align*}
\\
Also,
\begin{align*}	
	& \prod_{k=1}^{\alpha N \delta_{N}} P\left( \frac{1}{\sqrt{\alpha N}}\sum_{j=1}^{\alpha N + k}(Z_{j}-1) \geq - \phi_{N}(x) - \frac{k}{\sqrt{\alpha N}}  \right) \\
	& = \prod_{k=1}^{\alpha N \delta_{N}} P\left( \frac{1}{\sqrt{\alpha N + k}}\sum_{j=1}^{\alpha N + k}(Z_{j}-1) \geq - \left( \sqrt{\frac{\alpha N}{\alpha N + k}} \right)\left( \phi_{N}(x) + \frac{k}{\sqrt{\alpha N}} \right)  \right) \\
	& =  \prod_{k=1}^{\alpha N \delta_{N}} p_{k}
\end{align*}
where 
\begin{equation*}
	p_{k} = P\left( \frac{1}{\sqrt{\alpha N + k}}\sum_{j=1}^{\alpha N + k}(Z_{j}-1) \geq - \left( \sqrt{\frac{\alpha N}{\alpha N + k}} \right)\left( \phi_{N}(x) + \frac{k}{\sqrt{\alpha N}} \right)  \right)
\end{equation*}\\
\\
and, for fixed $k \in \left[ 0, \sqrt{2\alpha N \log(N)} \right]$, the factor $ \sqrt{\frac{\alpha N}{\alpha N + k}}$ goes to one as $N$ goes to infinity.\\ 
\\
\text{}\\
Finally,
\begin{align*}
	& \lim_{N \rightarrow +\infty}\log \left(\mathcal{P}_{N}\left(1 - \frac{f_{N}(x)}{\sqrt{\alpha N}}\right) \right) \\
	& = \lim_{N \rightarrow +\infty} \sum_{k=1}^{\sqrt{2\alpha N \log{N}}}   \log P\left( \frac{1}{\sqrt{\alpha N+k}}\sum_{j=1}^{\alpha N+k}(Z_{j}-1) \geq -\phi_{N}(x) - \frac{k}{\sqrt{\alpha N}}  \right)
\end{align*}
where $\delta_{N} = \sqrt{\frac{2\log{N}}{\alpha N}}$.\\ \\
\\
The Edgeworth expansion is then used to get the probability density of the standardised random variable  
\begin{align*}
	\frac{\sqrt{\alpha N + k}(\bar{Z}- \mu)}{\sigma} 
		& = \frac{1}{\sqrt{\alpha N+k}}\sum_{j=1}^{\alpha N+k}(Z_{j}-1)
\end{align*} 
where each element of sequence $Z_{1}, \cdots, Z_{\alpha N+k}$ is i.i.d. exponentially distributed with parameter equal to $1$.\text{}\\
\text{}\\
The variable $\bar{Z} = \frac{1}{\alpha N+k}\sum_{j=1}^{\alpha N+k}Z_{j}$ is the empirical mean. The mean of $\bar{Z}$ is equal to $1$ and its standard deviation is $\frac{\sigma}{\sqrt{\alpha N+k}} = \frac{1}{\sqrt{\alpha N+k}}$ with $\sigma = 1$.
\text{}\\
\text{}\\
Applying the Edgeworth expansion, as in \cite{Rider2003}, the logarithm of the probability $p_{k}$ is 
\begin{align*}
	& \log \left[  P\left( \frac{1}{\sqrt{\alpha N+k}}\sum_{j=1}^{\alpha N+k}(Z_{j}-1) \geq - \phi_{N}(x) - \frac{k}{\sqrt{\alpha N}}  \right) \right] \\
	& = \log\left[\int_{-\phi_{N}(x) - \frac{k}{\sqrt{\alpha N}} }^{+K_{N}} \frac{e^{-\frac{t^{2}}{2}}}{\sqrt{2\pi}}dt\right] + O\left( \frac{1}{\sqrt{\alpha N}}\sup_{\vert c\vert \leq H}\phi_{N}^{2}(c)e^{-\frac{\phi_{N}^{2}(c)}{2}}\right) + O\left(\frac{1}{\alpha N }\right)\\
	& + O\left(K_{N}(\alpha N)^{-3/2}\right) + O\left( e^{-\frac{1}{2}K_{N}^{2}} \right)
\end{align*}
where, similarly stated in \cite{Rider2003}, for $x$ restricted as in $\vert x \vert \leq H $ for some large positive $H$, $K_{N} $ which goes to $ + \infty$ faster than $\sup_{ \vert x \vert \leq H}\phi_{N}(x)$.\\ 
\\
Defining the function $f(k)$ as
\begin{equation*}
	f(k) = \log \left[ \int_{-\phi_{N}(x) - \frac{k}{\sqrt{\alpha N}}}^{+\infty} \frac{e^{-\frac{s^{2}}{2}}}{\sqrt{2\pi}}ds \right]
\end{equation*}
\\
Then,
\begin{equation*}
	\sum_{k=1}^{\sqrt{2\alpha N \log N}} f(k) = \int_{1}^{\sqrt{2\alpha N \log N}}\log \left[ \int_{-\phi_{N}(x) -\frac{k}{\sqrt{\alpha N}}}^{+\infty} \frac{e^{-\frac{s^{2}}{2}}}{\sqrt{2\pi}} \right]dt + E_{N}
\end{equation*}\\ \\
Let $u = \phi_{N}(x) + \frac{t}{\sqrt{\alpha N}}$ which implies that $dt = \sqrt{\alpha N} du$.\\
\begin{equation*}
	\sum_{k=1}^{\sqrt{2\alpha N \log N}} f(k) =  \sqrt{\alpha N} \int_{\phi_{N}(x) + \frac{1}{\sqrt{\alpha N}}}^{\phi_{N}(x) + \sqrt{2 \log N}}\log \left[ \int_{-u}^{+\infty} \frac{e^{-\frac{s^{2}}{2}}}{\sqrt{2\pi}}ds   \right]du + E_{N}
\end{equation*}
\\
where, with respect to the error terms $O\left( \frac{1}{\sqrt{\alpha N}}\sup_{\vert c\vert \leq H}\phi_{N}^{2}(c)e^{-\frac{\phi_{N}^{2}(c)}{2}}\right)$ and $ O\left(\frac{1}{\alpha N }\right)$ of the Edgeworth expansion of the logarithm of the probability $p_{k}$,
\begin{equation*}
	E_{N} = O\left(\left(\frac{\sqrt{\log N}}{\sqrt{\alpha N}} \right) \vee  \left( \sqrt{\log{N}}\sup_{\vert c \vert \leq H}\phi_{N}^{2}(c)e^{-\frac{\phi_{N}^{2}(c)}{2}} \right) \right)
\end{equation*}\\
\\
The term $\frac{1}{\sqrt{\alpha N}}$ appearing in the lower limit of integration is negligible. Extending the upper limit to $+\infty$, the function $\phi_{N}(x) = o \left( \sqrt{\log \left( \alpha N \right)}\right)$. This implies the following result
\begin{equation*}
	\lim_{N \rightarrow +\infty}\log \mathcal{P}_{N}\left( 1- \frac{f_{N}(x)}{\sqrt{\alpha N}} \right)  = \sqrt{\alpha N}\int_{\phi_{N}(x)}^{+\infty}\log\left[ \int_{-t}^{+\infty} \frac{e^{-\frac{s^{2}}{2}}}{\sqrt{2\pi}}ds\right]dt  + E_{N}
\end{equation*}
\\
Let $F_{\infty}(x)$ denote the limiting cumulative distribution function of the scaled minimum modulus of an $N \times N$ complex induced Ginibre matrix $G$ as $N$ goes to infinity. Then,
\begin{equation*}
	\log\left[ \int_{-t}^{+\infty} \frac{e^{-\frac{s^{2}}{2}}}{\sqrt{2\pi}}ds \right] = \log \left[ 1 - \int_{-\infty}^{-t} \frac{e^{-\frac{s^{2}}{2}}}{\sqrt{2\pi}}ds  \right]
\end{equation*}
and $\lim_{N \rightarrow +\infty} E_{N} = 0$.\\
\\
Thus, as the function $f_{N}$ goes to infinity, the logarithm of the limiting survival distribution function  $1-F_{\infty}(x)$ is 
\begin{align*}
		\log \left( 1 - F_{\infty}(x) \right)
		& = - \lim_{N \rightarrow +\infty} \sqrt{\frac{\alpha N}{2 \pi}} \int_{\phi_{N}(x)}^{+\infty}\left[ \int_{-\infty}^{-t} e^{-\frac{s^{2}}{2}}ds + O\left( \left( \int_{-\infty}^{-t} e^{-\frac{s^{2}}{2}} ds\right)^{2}\right) \right]dt \\
		& = - \lim_{N \rightarrow +\infty} \sqrt{\frac{\alpha N}{2 \pi}} \int_{\phi_{N}(x)}^{+\infty}\left[ \int_{t}^{+\infty} e^{-\frac{s^{2}}{2}}ds + O\left(\left( \int_{t}^{+\infty}e^{-\frac{s^{2}}{2}} ds\right)^{2} \right) \right]dt \\
	 	& = -\lim_{N \rightarrow +\infty} \sqrt{\frac{\alpha N}{2\pi}} \int_{\phi_{N}(x)}^{+\infty}\left[
	 	\frac{e^{-t^{2}/2}}{t} + O\left( \frac{e^{-t^{2}}}{t^{2}} \right)	\right]dt 	 	
\end{align*}\\
\\
Furthermore, $\phi_{N}(x) = 2f_{N}(x) - \frac{f_{N}^{2}(x)}{\sqrt{\alpha N}} = 2f_{N}(x)\left( 1 - \frac{f_{N}(x)}{2\sqrt{\alpha N}} \right)$ and the function $f_{N}(x) = o(\sqrt{\alpha N})$. This implies that 
\begin{align*}
	\lim_{N \rightarrow +\infty} \phi_{N}^{2}(x) 
	& = \lim_{N \rightarrow +\infty} 4f_{N}^{2}(x)\left( 1 - \frac{f_{N}(x)}{2\sqrt{\alpha N}} \right)^{2}\\
	& = \lim_{N \rightarrow +\infty} 4f_{N}^{2}(x)\left( 1 - \frac{f_{N}(x)}{\sqrt{\alpha N}} + \frac{f_{N}^{2}(x)}{4\alpha N} \right)
	  = \lim_{N \rightarrow +\infty} 4f_{N}^{2}(x)
\end{align*}
\\
Also,
\begin{equation*}
	\int_{\phi_{N}(x)}^{+\infty}
	 	\frac{1}{t}e^{-t^{2}/2} dt  =  \frac{1}{\phi_{N}^{2}(x)}e^{-\frac{\phi_{N}^{2}(x)}{2}} - \int_{\phi_{N}(x)}^{+\infty}\frac{2e^{-t^{2}/2}}{t^{3}}dt
\end{equation*}
\\
and in the limit as $N$ goes to infinity,
\begin{equation*}
	\int_{\phi_{N}(x)}^{+\infty}\frac{e^{-t^{2}/2}}{t^{3}}dt = \frac{1}{\phi_{N}^{2}(x)}e^{-\frac{\phi_{N}^{2}(x)}{2}} \times \frac{1}{\phi_N^{2}(x)} + O\left( \int_{\phi_{N}(x)}^{+\infty} 	
	\frac{e^{-t^{2}/2}}{t^{5}}dt \right)
\end{equation*}
This implies that
\begin{align*}
	& \lim_{N \rightarrow +\infty}\log \left(\mathcal{P}_{N}\left(1 - \frac{f_{N}(x)}{\sqrt{\alpha N}}\right) \right) 
	  = \log \left( 1 - F_{\infty}(x) \right) \\
	& = -\lim_{N \rightarrow +\infty} \sqrt{\frac{\alpha N}{2\pi}} \left[
	 	\frac{1}{\phi_{N}^{2}(x)}e^{-\frac{\phi_{N}^{2}(x)}{2}} \left(1 + O\left(\frac{1}{\phi_N^{2}(x)}\right)\right) \right] \\
	& = -\lim_{N \rightarrow +\infty} \sqrt{\frac{\alpha N}{2\pi}} \left[
	 	\frac{1}{4f_{N}^{2}(x)}e^{-2f_{N}^{2}(x)}  \left(1 + O\left(\frac{1}{4f_N^{2}(x)}\right)\right) \right]
\end{align*}\\
\\
For a convenient choice of the function $f_{N}(x)$, the limiting cumulative distribution function $F_{\infty}(x)$ is from  the class of Extreme Value Distributions composed of the three types of extreme value distributions for minima \cite{Gnedenko1943}.\\ \\ More precisely, choosing
\begin{equation*}
	f_{N}^{2}(x) = \frac{1}{2}\log \left( \frac{e^{-x}\sqrt{\alpha N / 2\pi}}{\log N} \right)
\end{equation*}
this implies
\begin{align*}
	& \lim_{N \rightarrow +\infty} \log \mathcal{P}_{N}\left[ r_{N} \geq 1 - \sqrt{\frac{1}{2 \alpha N}}\left( \log \frac{\sqrt{\alpha N/2\pi}}{\log{N}} - x \right)^{1/2}\right] \\
	& = -\lim_{N \rightarrow +\infty} \sqrt{\frac{\alpha N}{2\pi}} \left[
	 	\frac{1}{4f_{N}^{2}(x)}e^{-2f_{N}^{2}(x)} \left(1 + O\left(\frac{1}{4f_N^{2}(x)}\right)\right) \right]\\
	& = -\lim_{N \rightarrow +\infty}\exp(x)\frac{\log{N}}{2\log(e^{-x}\sqrt{\alpha N/2\pi} \times \frac{1}{\log{N}})}
\end{align*}\text{}\\
\text{}\\
Furthermore, for $\vert x \vert \ll N$ and using the same reasoning applied for the derivation of the limiting cumulative distribution function of the scaled spectral radius for matrices from the complex induced Ginibre ensemble
\begin{equation*}
	\lim_{N \rightarrow +\infty}\frac{\log{N}}{2\log(e^{-x}\sqrt{\alpha N/2\pi} \times \frac{1}{\log{N}})} = 1
\end{equation*}
\\
Finally,
\begin{equation*}
	\lim_{N \rightarrow +\infty} \log \mathcal{P}_{N}\left[ r_{N} \geq 1 - \sqrt{\frac{1}{2\alpha N}}\left( \log \frac{\sqrt{\alpha N/2\pi}}{\log{N}} - x \right)^{1/2}\right] = -\exp(x)
\end{equation*}\\
which is the logarithm of the survival distribution function of the standard Gumbel distribution for minima.\text{}\\
\text{}\\
Furthermore, setting $\gamma_{\alpha, N} =  \log \frac{\sqrt{\alpha N/2\pi}}{\log{N}} = \log{\sqrt{\alpha N/2\pi}} - \log\log{N}$ and using the Taylor expansion of the square root function\\
\begin{equation*}
	\left( \log \frac{\sqrt{\alpha N/2\pi}}{\log{N}} - x \right)^{1/2} 
	=  \gamma_{\alpha, N}^{1/2} \left( 1 - \frac{x}{\gamma_{\alpha, N}} \right)^{1/2}
	= \gamma_{\alpha, N}^{1/2} - \frac{x}{2\gamma_{\alpha, N}^{1/2}} + O\left(\left(\frac{x}{\gamma_{\alpha, N}}\right)^{2}\right)
\end{equation*}
\\
This implies that
\begin{align*}
	&\lim_{N \rightarrow +\infty} \mathcal{P}_{N}\left[ r_{N} \geq 1  - \left(\sqrt{\frac{\gamma_{\alpha, N}}{2\alpha N}} -\frac{x}{2\sqrt{2 \alpha N\gamma_{\alpha, N}}} \right) \right] 
	  = \exp(-\exp(x))
\end{align*}\\
\\
This limit corresponds to the survival distribution function of the standard Gumbel distribution for minima.\text{}\\
\text{}\\
The scaled minimum modulus $r_{N} = \frac{r_{min}^{(N)}(G)}{\sqrt{\alpha N}}$ is a random variable well approximated by
\begin{align*}
	r_{N} & \simeq  1 - \sqrt{\frac{\gamma_{\alpha, N}}{2 \alpha N}} + \frac{X}{2\sqrt{2 \alpha N\gamma_{\alpha, N}}} \\
		  & = 1 - \sqrt{\frac{\gamma_{\alpha, N}}{2 \alpha N}} - \frac{\log(Z)}{2\sqrt{2 \alpha N\gamma_{\alpha, N}}} \\
		  & = 1 - T_{\alpha, N} + \xi_{\alpha, N}
\end{align*}\\
\\
where the term $ T_{\alpha, N} = \sqrt{\frac{\gamma_{\alpha, N}}{2 \alpha N}}$ with $\gamma_{\alpha N} = \log{\sqrt{\alpha N/2\pi}} - \log\log{N}$. The random variable $Z$ denotes a random variable following a standard Exponential distribution which implies that $X = -\log(Z)$ has a standard Gumbel distribution. The random variable $\xi_{\alpha, N} = -\frac{\log(Z)}{2\sqrt{2 \alpha N\gamma_{\alpha, N}}}$ has a Gumbel distribution. \\ \\
\\ 
This completes the proof of Theorem \ref{RMinLimitTheoremComplexInducedGinibre}.
\\
\subsection{Proof of Theorem 10}\label{IndependenceRMaxRminComplexInducedGinibre}\text{}\\
\text{}\\
Let $G$ denote an $N \times N $ complex induced Ginibre matrix. Let $r_{in}= \sqrt{L}$ denote the inner radius and $r_{out}= \sqrt{L+N}$ the outer radius defining the edge of the eigenvalues support of the matrix $G$. The rectangularity index $L$ is proportional to $N$, i.e., $L = \alpha N$, $\alpha > 0$. Considering the scaling $\sqrt{\alpha N}$, the inner radius is now $r_{in}= 1$ and the outer radius $r_{out}= \sqrt{\frac{1+\alpha}{\alpha}} = \rho$.\\
\\
The independence of the scaled spectral radius $ R_{N} = \frac{r_{max}^{(N)}(G)}{\sqrt{\alpha N}}$ and the scaled minimum modulus $r_{N} =\frac{r_{min}^{(N)}(G)}{\sqrt{\alpha N}}$ is derived from their joint probability 
\begin{align*}
P\left(r_{N} \geq r  \text{}\text{ and }\text{} R_{N} \leq R \right) 
	 & = P\left(r_{min}^{(N)}(G)\geq \sqrt{\alpha N }r \text{}\text{ and }\text{} r_{max}^{(N)}(G) \leq \sqrt{\alpha N }R\right) \\
	 & = \frac{1}{\prod_{k=1}^{N}\Gamma(k + \alpha N)} \prod_{k=0}^{N-1} \int_{\alpha Nr^{2}}^{\alpha NR^{2}} t^{k+L} e^{-t} dt \\
 	 & = \prod_{k=1}^{N}\int_{\alpha Nr^{2}}^{\alpha NR^{2}}f_{Gamma(k+L, 1)}(t)dt
\end{align*}
\\
The function $f_{Gamma(k+L, 1)}(t)$ is the probability density function of the Gamma distribution with shape parameter $k+L$ and rate parameter $1$.\\
\\
Thus,
\begin{align*}
	P\left(r_{min}^{(N)}(G) \geq \sqrt{\alpha N }r \text{}\text{ and }\text{} r_{max}^{(N)}(G) \geq \sqrt{\alpha N }R\right)
   	 & = \prod_{k=1}^{N} P\left( r^{2} \leq \frac{1}{\alpha N}\sum_{j=1}^{\alpha N +k} Z_{j} \leq R^{2} \right)
\end{align*} 
where the $Z_{j}$ for $ j = \lbrace  1, \cdots, k+L  \rbrace$ are independent and identically distributed random variables following a standard Exponential distribution.\\ 
\\
Let $r = 1 - \frac{f_{N}(y)}{\sqrt{\alpha N}}$ and $R =  \sqrt{\frac{1+\alpha}{\alpha}} + \frac{f_{N}(x)}{\sqrt{\alpha N}}= \rho + \frac{f_{N}(x)}{\sqrt{\alpha N}} $ where $f_{N}(x)$ and $f_{N}(y)$ are $o\left( \sqrt{\alpha N} \right)$.\text{}\\ \text{}\\ 
\text{}\\
This implies
\begin{align*}
	& P\left(r_{min}^{(N)}(G) \geq \sqrt{\alpha N }r \text{}\text{ and }\text{} r_{max}^{(N)}(G) \leq \sqrt{\alpha N }R\right)\\
   	& = \prod_{k=1}^{N} P\left( r^{2} \leq \frac{1}{\alpha N}\sum_{j=1}^{k+L}Z_{j} \leq R^{2} \right)\\
	& = \prod_{k=1}^{N} \left[P\left( -\phi_{N}^{(r)}(y) - \frac{k}{\sqrt{\alpha N}} \leq  \frac{1}{\sqrt{\alpha N}}\sum_{j=1}^{k+L}\left( Z_{j} -1 \right) \leq \phi_{N}^{(R)}(x) + C_{\alpha,N} - \frac{k}{\sqrt{\alpha N}} \right)\right]
\end{align*}\text{}\\
\text{}\\
where $C_{\alpha,N} = \sqrt{\alpha N}\left(\rho^{2}-1\right) > 0$. The functions $\phi_{N}^{(r)}(y) = 2f_{N}(y)\left( 1 - \frac{f_{N}(y)}{2\sqrt{\alpha N}} \right)$ and $\phi_{N}^{(R)}(x)  =  2\rho f_{N}(x)\left( 1 - \frac{f_{N}(x)}{2 \rho \sqrt{\alpha N}} \right)$.\\ \\ 
\\
Furthermore, with a proportional rectangularity index $L = \alpha N$, $\alpha > 0$, the partial product\\ 
\begin{align*}
	& \prod_{k=1}^{N} P\left( -\phi_{N}^{(r)}(y) - \frac{k}{\sqrt{\alpha N}} \leq  \frac{1}{\sqrt{\alpha N}}\sum_{j=1}^{\alpha N + k}\left( Z_{j} -1 \right) \leq \phi_{N}^{(R)}(x) + C_{\alpha,N} - \frac{k}{\sqrt{\alpha N}} \right)
\end{align*}\\
is bounded by any partial product of $p_{k}$, i.e., $\prod_{k=1}^{\alpha N \delta_{N}} p_{k}$, for a positive $\delta_{N}$ less than 1 (cf. \cite{Rider2003}), such as \\
\begin{equation*}
	\prod_{k=1}^{\alpha N \delta_{N}} p_{k} \prod_{k=\alpha N \delta_{N}}^{N} p_{k}
	\leq P\left( -\phi_{N}^{(r)}(y) - \frac{k}{\sqrt{\alpha N}} \leq \frac{\sum_{j=1}^{k+L}\left( Z_{j} -1 \right)}{\sqrt{\alpha N}} \leq \phi_{N}^{(R)}(x) + C_{\alpha,N} - \frac{k}{\sqrt{\alpha N}} \right)  \leq \prod_{k=1}^{\alpha N \delta_{N}} p_{k}
\end{equation*}\\
where	
\begin{equation*}
	p_{k} = P\left( -\phi_{N}^{(r)}(y) - \frac{k}{\sqrt{\alpha N}} \leq  \frac{1}{\sqrt{\alpha N}}\sum_{j=1}^{k+L}\left( Z_{j} -1 \right) \leq \phi_{N}^{(R)}(x) + C_{\alpha,N} - \frac{k}{\sqrt{\alpha N}} \right)
\end{equation*}\\ 
\\
Furthermore, \\
\small \begin{align*}
	 & \prod_{k=\alpha N \delta_{N}}^{N} P\left( -\phi_{N}^{(r)}(y) - \frac{k}{\sqrt{\alpha N}} \leq \frac{\sum_{j=1}^{k+L}\left( Z_{j} -1 \right)}{\sqrt{\alpha N}} \leq \phi_{N}^{(R)}(x) + C_{\alpha,N} - \frac{k}{\sqrt{\alpha N}} \right) \\
	 & = \prod_{k=\alpha N \delta_{N}}^{N} \left[ P\left( \frac{\sum_{j=1}^{k+L}\left( Z_{j} -1 \right)}{\sqrt{\alpha N}} \leq \phi_{N}^{(R)}(x) + C_{\alpha,N} - \frac{k}{\sqrt{\alpha N}} \right) 
	 - P\left(   \frac{\sum_{j=1}^{k+L}\left( Z_{j} - 1 \right)}{\sqrt{\alpha N}} \leq -\phi_{N}^{(r)}(y) - \frac{k}{\sqrt{\alpha N}} \right) \right] \\
	 & = \prod_{k=\alpha N \delta_{N}}^{N} \left[ 1 -  P\left( \frac{\sum_{j=1}^{k+L}\left( Z_{j} -1 \right)}{\sqrt{\alpha N}} \geq \Phi_{N}^{(R)}(x) - \frac{k}{\sqrt{\alpha N}} \right) 
	 - P\left(   \frac{\sum_{j=1}^{k+L}\left( Z_{j} - 1 \right)}{\sqrt{\alpha N}} \leq -\phi_{N}^{(r)}(y) - \frac{k}{\sqrt{\alpha N}} \right) \right] \\
	 & = \prod_{k=\alpha N \delta_{N}}^{N} \left[ 1 -  P\left( \sum_{j=1}^{k+L}\left( Z_{j} -1 \right)   \geq \sqrt{\alpha N}\Phi_{N}^{(R)}(x) - k \right) 
	 - P\left(   \frac{\sum_{j=1}^{k+L}\left( Z_{j} - 1 \right)}{\sqrt{\alpha N}} \leq -\phi_{N}^{(r)}(y) - \frac{k}{\sqrt{\alpha N}} \right) \right] \\
	 & = \prod_{k=\alpha N \delta_{N}}^{N} \left[ 1 -  P\left( \sum_{j=1}^{k+L} Z_{j}  \geq \sqrt{\alpha N}\Phi_{N}^{(R)}(x) + \alpha N \right) 
	 - P\left( \frac{\sum_{j=1}^{k+L} Z_{j}}{\sqrt{\alpha N + k}} \leq - \left( \sqrt{\frac{\alpha N}{\alpha N + k}} \right)\phi_{N}^{(r)}(y) +  \frac{k}{\alpha N + k} \right) \right] \\
	 & > \prod_{k=\alpha N \delta_{N}}^{N} \left[ 1 -  P\left( \sum_{j=1}^{k+L} Z_{j}  \geq \sqrt{\alpha N}\Phi_{N}^{(R)}(x) + \alpha N \right) 
	 - P\left( \frac{\sum_{j=1}^{k+L} Z_{j}}{\sqrt{\alpha N + k}} \leq   \frac{k}{\alpha N + k} \right) \right] \\
	 & = \prod_{k=\alpha N \delta_{N}}^{N} \left[ 1 -  P\left( \sum_{j=1}^{k+L} Z_{j}  \geq \sqrt{\alpha N}\Phi_{N}^{(R)}(x) + \alpha N \right) 
	 - P\left( \frac{\sum_{j=1}^{k+L} Z_{j}}{\sqrt{\alpha N + k}} \leq   \frac{k}{\alpha N + k} \right) \right] \\
	 & > \prod_{k=\alpha N \delta_{N}}^{N} \left[ 1 -  O\left(\frac{1}{N}\right) \right] 
\end{align*}\\
\\
Now,\\
\small \begin{align*}
	& \log  \prod_{k=1}^{\alpha N\delta_{N}} P\left( -\phi_{N}^{(r)}(y) - \frac{k}{\sqrt{\alpha N}} \leq \frac{1}{\sqrt{\alpha N}}\sum_{j=1}^{k+\alpha N}\left( Z_{j} -1 \right) \leq \phi_{N}^{(R)}(x) + C_{\alpha,N} - \frac{k}{\sqrt{\alpha N}} \right)\\ 
	&  = \sum_{k=1}^{\sqrt{2\alpha N \log N}}\log P\left( \alpha_{k}^{(N)}\left(-\phi_{N}^{(r)}(y) - \frac{k}{\sqrt{\alpha N}}\right) \leq  \frac{1}{\sqrt{\alpha N+k}}\sum_{j=1}^{k+ \alpha N}\left( Z_{j} -1 \right) \leq \alpha_{k}^{(N)} \left(\phi_{N}^{(R)}(x) + C_{\alpha,N} - \frac{k}{\sqrt{\alpha N}}\right) \right)
\end{align*}
\normalsize where, for fixed $k \in \left[ 0, \alpha N \delta_{N} \right]$, the factor $ \alpha_{k}^{(N)} = \sqrt{\frac{\alpha N}{\alpha N + k}}$ goes to one as $N$ goes to infinity.\\ 
\\
Furthermore, using the classical Edgeworth expansion, in line with equation $(9)$ in \cite{Rider2003}, this implies that
\begin{align*}
	& \log P\left( -\phi_{N}^{(r)}(y) - \frac{k}{\sqrt{\alpha N}} \leq  \frac{1}{\sqrt{\alpha N + k}}\sum_{j=1}^{k+L}\left( Z_{j} -1 \right) \leq  \phi_{N}^{(R)}(x) + C_{\alpha,N} - \frac{k}{\sqrt{\alpha N}} \right)\\
	& = \log \left[\int_{-\phi_{N}^{(r)}(y) - \frac{k}{\sqrt{\alpha N}}}^{\phi_{N}^{(R)}(x) + C_{\alpha,N} - \frac{k}{\sqrt{\alpha N}}} \frac{e^{-\frac{t^{2}}{2}}}{\sqrt{2\pi}}  dt \right]
	  + O\left(\frac{1}{\sqrt{\alpha N}}\sup_{\vert c \vert  \leq  H} ( \phi_{N}^{(R)}(c))^{2}e^{-\frac{(\phi_{N}^{(R)}(c))^{2}}{2}}\right)\\
	& + O\left(\frac{1}{\alpha N}\right) + O\left( \phi_{N}^{(r)}(y)(\alpha N)^{-3/2}\right) + O\left(e^{-\frac{(\phi_{N}^{(r)}(y))^{2}}{2}}\right)
\end{align*}\\
where, as in \cite{Rider2003}, for $x$ restricted as in $ \vert x \vert \leq H$ for some large positive $H$, $K_{N}$ goes to $+\infty$ faster than $\sup_{\vert x \vert \leq H}\phi_{N}^{(R)}(x)$.\\
\\ 
\\
Then,
\begin{align*}
	& \sum_{k=1}^{\sqrt{2\alpha N\log N}}  \log \left[\int_{-\phi_{N}^{(r)}(y) - \frac{k}{\sqrt{\alpha N}}}^{\phi_{N}^{(R)}(x) + C_{\alpha,N} - \frac{k}{\sqrt{\alpha N}}} \frac{e^{-\frac{t^{2}}{2}}}{\sqrt{2\pi}}  dt \right]\\ 
	& = \sum_{k=1}^{\sqrt{2\alpha N\log N}}  \log \left[ 1 - \int_{\phi_{N}^{(R)}(x) + C_{\alpha,N} - \frac{k}{\sqrt{\alpha N}}}^{+\infty} \frac{e^{-\frac{t^{2}}{2}}}{\sqrt{2\pi}} dt  - \int_{-\infty}^{-\phi_{N}^{(r)}(y) - \frac{k}{\sqrt{\alpha N}}} \frac{e^{-\frac{t^{2}}{2}}}{\sqrt{2\pi}} dt  \right] 
\end{align*}\\
In the limit as $N$ goes to infinity, 
\begin{align*}
	& \sum_{k=1}^{\sqrt{2\alpha N\log N}}  \log \left[ 1 - \int_{\phi_{N}^{(R)}(x) + C_{\alpha,N} - \frac{k}{\sqrt{\alpha N}}}^{+\infty} \frac{e^{-\frac{t^{2}}{2}}}{\sqrt{2\pi}} dt  - \int_{-\infty}^{-\phi_{N}^{(r)}(y) - \frac{k}{\sqrt{\alpha N}}} \frac{e^{-\frac{t^{2}}{2}}}{\sqrt{2\pi}} dt  \right] \\
	& = \sum_{k=1}^{\sqrt{2\alpha N\log N}}\left[ - \int_{\phi_{N}^{(R)}(x) + C_{\alpha,N} - \frac{k}{\sqrt{\alpha N}}}^{+\infty} \frac{e^{-\frac{t^{2}}{2}}}{\sqrt{2\pi}} dt  - \int_{-\infty}^{-\phi_{N}^{(r)}(y) - \frac{k}{\sqrt{\alpha N}}} \frac{e^{-\frac{t^{2}}{2}}}{\sqrt{2\pi}} dt \right]\\
	& + \sum_{k=1}^{\sqrt{2\alpha N\log N}} O\left( \left(   \int_{\phi_{N}^{(R)}(x) + C_{\alpha,N} + \frac{k}{\sqrt{\alpha N}}}^{+\infty} \frac{e^{-\frac{t^{2}}{2}}}{\sqrt{2\pi}} dt  + \int_{-\infty}^{-\phi_{N}^{(r)}(y) - \frac{k}{\sqrt{\alpha N}}} \frac{e^{-\frac{t^{2}}{2}}}{\sqrt{2\pi}} dt  \right)^{2} \right) \\
	& = -\sum_{k=1}^{\sqrt{2\alpha N\log N}}  \int_{\Phi_{N}^{(R)}(x) - \frac{k}{\sqrt{\alpha N}}}^{+\infty} \frac{e^{-\frac{t^{2}}{2}}}{\sqrt{2\pi}} dt 
	  \times \left( 1 + O\left(\int_{\Phi_{N}^{(R)}(x) - \frac{k}{\sqrt{\alpha N}}}^{+\infty} \frac{e^{-\frac{t^{2}}{2}}}{\sqrt{2\pi}} dt \right) + O\left(\int_{-\infty}^{-\phi_{N}^{(r)}(y) - \frac{k}{\sqrt{\alpha N}}} \frac{e^{-\frac{s^{2}}{2}}}{\sqrt{2\pi}} ds \right) \right)\\
	&  - \sum_{k=1}^{\sqrt{2\alpha N\log N}}  \int_{-\infty}^{-\phi_{N}^{(r)}(y) - \frac{k}{\sqrt{\alpha N}}} \frac{e^{-\frac{t^{2}}{2}}}{\sqrt{2\pi}} dt 
	  \times \left( 1 + O\left(\int_{-\infty}^{-\phi_{N}^{(r)}(y)  - \frac{k}{\sqrt{\alpha N}}} \frac{e^{-\frac{t^{2}}{2}}}{\sqrt{2\pi}} dt \right)  \right)\\
	& = -\sum_{k=1}^{\sqrt{2\alpha N\log N}}  \int_{\Phi_{N}^{(R)}(x) - \frac{k}{\sqrt{\alpha N}}}^{+\infty} \frac{e^{-\frac{t^{2}}{2}}}{\sqrt{2\pi}} dt 
	  \times \left( 1 + O\left(\int_{\Phi_{N}^{(R)}(x) - \frac{k}{\sqrt{\alpha N}}}^{+\infty} \frac{e^{-\frac{t^{2}}{2}}}{\sqrt{2\pi}} dt \right) + O\left( \frac{e^{-\frac{\left( \phi_{N}^{(r)}(y) \right)^{2}}{2}}}{\phi_{N}^{(r)}(y)} \right) \right)\\
	&  - \sum_{k=1}^{\sqrt{2\alpha N\log N}}  \int_{-\infty}^{-\phi_{N}^{(r)}(y) - \frac{k}{\sqrt{\alpha N}}} \frac{e^{-\frac{t^{2}}{2}}}{\sqrt{2\pi}} dt 
	  \times \left( 1  + O\left(\int_{-\infty}^{-\phi_{N}^{(r)}(y) - \frac{k}{\sqrt{\alpha N}}} \frac{e^{-\frac{t^{2}}{2}}}{\sqrt{2\pi}} dt \right)  \right)
\end{align*}\\
This implies that,
\begin{align*}
	& \sum_{k=1}^{\sqrt{2\alpha N\log N}}  \log \left[ 1 - \int_{\phi_{N}^{(R)}(x) + C_{\alpha,N} - \frac{k}{\sqrt{\alpha N}}}^{+\infty} \frac{e^{-\frac{t^{2}}{2}}}{\sqrt{2\pi}} dt  - \int_{-\infty}^{-\phi_{N}^{(r)}(y) - \frac{k}{\sqrt{\alpha N}}} \frac{e^{-\frac{t^{2}}{2}}}{\sqrt{2\pi}} dt  \right] \\
	& = -\sum_{k=1}^{\sqrt{2\alpha N\log N}}  \int_{\Phi_{N}^{(R)}(x) - \frac{k}{\sqrt{\alpha N}}}^{+\infty} \frac{e^{-\frac{t^{2}}{2}}}{\sqrt{2\pi}} dt 
	  \times \left( 1 + O\left(\int_{\Phi_{N}^{(R)}(x) - \frac{k}{\sqrt{\alpha N}}}^{+\infty} \frac{e^{-\frac{t^{2}}{2}}}{\sqrt{2\pi}} dt \right) \right)\\
	& - \sum_{k=1}^{\sqrt{2\alpha N\log N}}  \int_{-\infty}^{-\phi_{N}^{(r)}(y) - \frac{k}{\sqrt{\alpha N}}} \frac{e^{-\frac{t^{2}}{2}}}{\sqrt{2\pi}} dt 
	  \times \left( 1 + O\left( \int_{-\infty}^{-\phi_{N}^{(r)}(y) - \frac{k}{\sqrt{\alpha N}}} \frac{e^{-\frac{t^{2}}{2}}}{\sqrt{2\pi}} dt \right)  \right)
\end{align*}\\
where $\Phi_{N}^{(R)}(x) = \phi_{N}^{(R)}(x) + C_{\alpha,N}$.\\
\\
\begin{remark} For $\alpha = O\left(N^{m} \right), (\rho^{2} - 1)\sqrt{\alpha N} = \frac{1}{\alpha}\sqrt{\alpha N} = O \left( N^{\frac{1-m}{2}} \right)$. There exists a constant $ m_{1} \geq 1$, $\forall m > m_{1}$ and $\alpha = O\left( N^{m}\right)$, as $N$ goes to infinity, $\phi_{N}^{(R)}(x) = o\left( \sqrt{\alpha N} \right) > \sqrt{2\log(N)}$. This assumption is used such that the term $C_{\alpha,N}$, goes to zero as $N$ goes to infinity.
\end{remark}\text{}\\
\\
The leading order sum as the Riemann integral of last equation is,\\
\begin{align*}
	&-\sum_{k=1}^{\sqrt{2\alpha N\log N}}  \int_{\Phi_{N}^{(R)}(x) - \frac{k}{\sqrt{\alpha N}}}^{+\infty} \frac{e^{-\frac{t^{2}}{2}}}{\sqrt{2\pi}} dt 
	  \times \left(  1 + O\left( \int_{\Phi_{N}^{(R)}(x) - \frac{k}{\sqrt{\alpha N}}}^{+\infty} \frac{e^{-\frac{t^{2}}{2}}}{\sqrt{2\pi}} dt \right) + O\left( \frac{e^{-\frac{\left( \phi_{N}^{(r)}(y) \right)^{2}}{2}}}{\phi_{N}^{(r)}(y)} \right) \right)\\
	&  - \sum_{k=1}^{\sqrt{2\alpha N\log N}}  \int_{-\infty}^{-\phi_{N}^{(r)}(y) - \frac{k}{\sqrt{\alpha N}}} \frac{e^{-\frac{t^{2}}{2}}}{\sqrt{2\pi}} dt 
	  \times \left( 1  + O\left( \int_{-\infty}^{-\phi_{N}^{(r)}(y) - \frac{k}{\sqrt{\alpha N}}} \frac{e^{-\frac{t^{2}}{2}}}{\sqrt{2\pi}} dt \right)  \right)\\
	& = - \sqrt{\alpha N} \int_{\Phi_{N}^{(R)}(x)}^{\Phi_{N}^{(R)}(x) - \sqrt{2\log(N)}}\left( \int_{u}^{+\infty}\frac{e^{-\frac{t^{2}}{2}}}{\sqrt{2\pi}} dt \times \left( 1 +  O\left( \int_{u}^{+\infty}\frac{e^{-\frac{t^{2}}{2}}}{\sqrt{2\pi}} dt \right) \right)\right) du  \\
	& -\left(\sum_{k=1}^{\sqrt{2\alpha N\log N}}  \int_{\Phi_{N}^{(R)}(x) - \frac{k}{\sqrt{\alpha N}}}^{+\infty} \frac{e^{-\frac{t^{2}}{2}}}{\sqrt{2\pi}} dt \right) \times 	O\left( \frac{e^{-\frac{\left( \phi_{N}^{(r)}(y) \right)^{2}}{2}}}{\phi_{N}^{(r)}(y)} \right)\\
	& - \sqrt{\alpha N} \int_{\phi_{N}^{(r)}(y)}^{\phi_{N}^{(r)}(y) + \sqrt{2\log(N)}} \left(\int_{-\infty}^{-u}\frac{e^{-\frac{t^{2}}{2}}}{\sqrt{2\pi}} dt \times \left( 1 + O\left(\int_{-\infty}^{-u}\frac{e^{-\frac{t^{2}}{2}}}{\sqrt{2\pi}} dt \right) \right) \right) du  + E_{N}
\end{align*}
where 
\begin{equation*}
	E_{N} = O\left(  \left(\sqrt{\frac{\log N}{\alpha N}} \right) \vee  \left( \sqrt{\log{N}}\sup_{\vert c \vert \leq H}(\Phi_{N}^{(R)}(c))^{2}e^{-\frac{(\Phi_{N}^{(R)}(c))^{2}}{2}} \right) \right)
\end{equation*}\\
\\
\\
Also, $\lim_{N \longrightarrow +\infty} \frac{e^{-\frac{\left( \phi_{N}^{(r)}(y) \right)^{2}}{2}}}{\phi_{N}^{(r)}(y)} = 0$. Thus, following the framework presented in \cite{Rider2003}, choosing,
\begin{equation*}
	f_{N}^{2}(x) = \frac{1}{2\rho^{2}}\left( \log\left( \frac{e^{x}\sqrt{N/(2\pi)}}{\log(N)} \right) \right)
\end{equation*}
\\
this implies\\
\begin{align*}
	\log F_{X}^{(\infty)}(x) 
		& = - \lim_{N \rightarrow +\infty} \sqrt{\alpha N} \int_{\phi_{N}^{(R)}(x)}^{+\infty}\left( \int_{u}^{+\infty}\frac{e^{-\frac{t^{2}}{2}}}{\sqrt{2\pi}} dt \times \left( 1 + O\left(\int_{u}^{+\infty}\frac{e^{-\frac{t^{2}}{2}}}{\sqrt{2\pi}} dt \right) \right)\right) du \\
		& =  - \lim_{N \rightarrow +\infty} \sqrt{\alpha N} \int_{\phi_{N}^{(R)}(x)}^{+\infty}\left( \int_{u}^{+\infty}\frac{e^{-\frac{t^{2}}{2}}}{\sqrt{2\pi}} dt \right) du  \times \left( 1 + O\left(\frac{1}{\left(\phi_{N}^{(R)}(x) \right)^{2}} \right) \right) \\
		& =  - \lim_{N \rightarrow +\infty} \sqrt{\frac{\alpha N}{2\pi}}\frac{e^{-\frac{\left(\phi_{N}^{(R)}(x)\right)^{2}}{2}}}{\left( \phi_{N}^{(R)}(x) \right)^{2}}  \times \left( 1 + O\left(\frac{1}{\left(\phi_{N}^{(R)}(x) \right)^{2}} \right) \right)  \\
		& =  -\exp(-x)
\end{align*}\\
\\
And, choosing 
\begin{equation*}
	f_{N}^{2}(y) = \frac{1}{2}\left( \log\left( \frac{e^{-y}\sqrt{N/(2\pi)}}{\log(N)} \right) \right)
\end{equation*}
\begin{align*}
	\log\left( 1- F_{Y}^{(\infty)}(y)\right) 
		& = - \lim_{N \rightarrow +\infty} \sqrt{\alpha N} \int_{\phi_{N}^{(r)}(y)}^{+\infty} \left(\int_{-\infty}^{-u}\frac{e^{-\frac{t^{2}}{2}}}{\sqrt{2\pi}} dt \times \left( 1 + O\left( \int_{-\infty}^{-u}\frac{e^{-\frac{t^{2}}{2}}}{\sqrt{2\pi}} dt \right) \right) \right) du \\
		& = -\exp(y)
\end{align*} 
this implies 
\begin{align*}
	& \lim_{N \rightarrow +\infty} P\left[ r_{N} \geq 1  - \sqrt{\frac{\gamma_{\alpha, N}}{2\alpha N}} + \frac{y}{2\sqrt{2 \alpha N\gamma_{\alpha, N}}}
	 \text{}\text{ and }\text{} R_{N} \leq  \rho  + \sqrt{\frac{\gamma_{\alpha, N}}{2\rho^{2}\alpha N}} +\frac{x}{2\sqrt{2\rho^{2}\alpha N\gamma_{\alpha, N}}}\right] \\
	& =  F_{X}^{(\infty)}(x) \left(  1 - F_{Y}^{(\infty)}(y)\right)
\end{align*}
with $\gamma_{\alpha, N} =  \log \frac{\sqrt{\alpha N/2\pi}}{\log{N}} = \log{\sqrt{\alpha N/2\pi}} - \log\log{N}$.\\ 
\\
\\
These limiting distribution functions are the standard Gumbel cumulative distribution function for maxima and the standard Gumbel survival distribution function for minima.\text{}\\
\text{}\\
This proves the independence of the scaled spectral radius $R_{N}=\frac{ r_{max}^{(N)}(G)}{\sqrt{\alpha N}}$ and the scaled minimum modulus $ r_{N}=\frac{r_{min}^{(N)}(G)}{\sqrt{\alpha N}}$ of an $N \times N$ matrix $G$ from the complex induced Ginibre ensemble with proportional rectangularity index $L$ as $N$ goes to infinity. \\ \\

\section{Conditional probability for the complex Ginibre ensemble}\label{ComplexGinibreConditionalProbability}\text{}\\
Let $\mathcal{D}(z_0, s)$ denote the disk of radius $s$ centred at $z_0$. The conditional probability $H^{(N)}(s, z_0)$, that given one eigenvalue lies at the point $z_0$ and the others are found outside $\mathcal{D}(z_0, s)$, is a result initially presented in \cite{GrobeHaakeSommers1988} and then referenced in \cite{KhoruzhenkoSommers2015}. This result is derived as follows. Let $z_{1}, ..., z_{N}$ denote  the $N$ eigenvalues of an $N \times N$ complex Ginibre matrix. The conditional probability $H^{(N)}(s, z_0)$ is then,\\
\begin{align*}
	H^{(N)}(s, z_0) 
	& = P(\left\lbrace z_2, \cdots, z_N \right\rbrace \notin \mathcal{D}(z_0, s) \vert z_1 = z_0 )\\
	& = \frac{P(\left\lbrace z_2, \cdots, z_N \right\rbrace \notin \mathcal{D}(z_0, s) \cap z_1 = z_0  )}{P(z_1 = z_0)}\\
	& = \frac{1}{P(z_1 = z_0)}\int_{\mathbb{C}}d^2z_2 \cdots \int_{\mathbb{C}}d^2z_NP(z_0, z_2, \cdots, z_N)\prod_{k=2}^{N}[1_{\left\lbrace z_k \notin \mathcal{D}(z_0, s)\right\rbrace }]
\end{align*}\\
The probability $P(z_1 = z_0)$ is the probability density at the point $z_0$. It is equal to the one-point correlation function in the vicinity of the point $z_0$ (i.e., the density in the vicinity of the point $z_0$) divided by the total number of eigenvalues which is equal to $N$.\text{}\\ \text{}\\
More precisely,
\begin{equation*}
	P(z_1 = z_0) = \frac{R_1^{(N)}(z_0)}{N}
\end{equation*}
This implies that
\begin{equation*}
	H^{(N)}(s, z_0) = \frac{N}{R_1^{(N)}(z_0)}\int_{\mathbb{C}}d^2z_2 \cdots \int_{\mathbb{C}}d^2z_NP(z_0, z_2, \cdots, z_N)\prod_{k=2}^{N}[1_{\left\lbrace z_k \notin \mathcal{D}(z_0, s)\right\rbrace }]
\end{equation*}
\text{}\\
The joint probability density function of the eigenvalues is \\
\begin{equation*}
	P(z_{1}, z_{2}, \cdots, z_{N}) = \frac{1}{N!\pi^{N}\prod_{j=0}^{N-1}j!}e^{-\sum_{j=1}^{N}\vert z_j \vert^2}\prod_{1 \leq i<j \leq N}\vert z_i - z_j \vert^2
\end{equation*}
\\
Let $z_1 = z_0 = 0$. Then,\\
\begin{equation*}
	\prod_{1\leq i < j \leq N}\vert z_i - z_j \vert^2= \prod_{j=2}^{N}\vert z_j \vert^2 \prod_{2 \leq i < j \leq N}\vert z_i - z_j \vert^2
\end{equation*}
\text{}\\
The term \\
\begin{align*}
	\prod_{2 \leq i < j \leq N}\vert z_i - z_j \vert^2  
		& = \vert \Delta(z_2, \cdots, z_N)\vert^2 \\
		& = \Delta(z_2, \cdots, z_N)\overline{ \Delta}(z_2, \cdots, z_N) 
		  = \det\left( z_{k}^{N-j}\right)_{k,j=2}^{N}\det\left( \bar{z}_{k}^{N-j}\right)_{k,j=2}^{N}
\end{align*}\\
which implies
\begin{align*}
	\prod_{j=2}^{N} \vert z_j \vert^2\prod_{2 \leq i < j \leq N}\vert z_i - z_j \vert^2 
		& = \prod_{j=2}^{N} \vert z_j \vert^2 \det\left( z_{k}^{N-j}\right)_{k,j=2}^{N}\det\left( \bar{z}_{k}^{N-j}\right)_{k,j=2}^{N}\\
		& = \prod_{j=2}^{N}z_j  \det\left( z_{k}^{N-j}\right)_{k,j=2}^{N} \prod_{j=2}^{N}\bar{z}_{j} \det\left( \bar{z}_{k}^{N-j}\right)_{k,j=2}^{N}\\
		& = \det\left( z_{k}^{N+1-j}\right)_{k,j=2}^{N} 
	\det\left( \bar{z}_{k}^{N+1-j}\right)_{k,j=2}^{N}
\end{align*}
Finally,
\begin{align*}
	P(0, z_2, \cdots, z_N) 
		& = \frac{1}{N!\pi^{N}\prod_{j=0}^{N-1}j!}e^{-\sum_{j=2}^{N}\vert z_j \vert^2}\prod_{1 \leq i<j \leq N}\vert z_i - z_j \vert^2\\
		& = \frac{1}{N!\pi^{N}\prod_{j=0}^{N-1}j!}e^{-\sum_{j=2}^{N}\vert z_j \vert^2}\det\left( z_{k}^{N+1-j}\right)_{k,j=2}^{N} \det\left( \bar{z}_{k}^{N+1-j}\right)_{k,j=2}^{N}
\end{align*}
\\
Furthermore,
\begin{equation*}
	R_1^{(N)}(0) = \frac{1}{\pi}\frac{\Gamma(N, 0)}{\Gamma(N)} = \frac{1}{\pi}
\end{equation*}
\text{}\\
\text{}\\
Applying Andreief's integration formula,
\begin{align*}
	H^{(N)}(s, 0) 
	& = \frac{N(N-1)!}{N!\pi^{N-1}\prod_{j=0}^{N-1}j!}\prod_{j=2}^{N}\left[\int_{\vert z \vert > s} \left( \vert z\vert^2 \right)^{N+1-j} e^{-\vert z\vert^2}d^2z\right]\\
	& = \frac{1}{\prod_{j=0}^{N-1}j!} \prod_{k=1}^{N-1} \int_{s^2}^{+ \infty} e^{- t}t^k dt
	  = \prod_{k=1}^{N-1}\frac{\Gamma(k+1, s^2)}{\Gamma(k+1)}		
\end{align*}\\
Also,
\begin{equation*}
	e^{x} = \sum_{k=0}^{+\infty}\frac{x^{k}}{k!} = \sum_{k=0}^{n}\frac{x^{k}}{k!} + \sum_{k=n+1}^{+\infty}\frac{x^{k}}{k!}
\end{equation*}
\\
and then,
\begin{equation*}
	e^{-x}\sum_{k=0}^{n}\frac{x^{k}}{k!} = 1 - e^{-x}\sum_{k=n+1}^{+\infty}\frac{x^{k}}{k!}
\end{equation*}
which implies that
\begin{equation*}
	H^{(N)}(s, 0) 
		= \prod_{k=1}^{N-1}\frac{\Gamma(k+1, s^2)}{\Gamma(k+1)}
		= \prod_{k=1}^{N-1}e^{-s^{2}}\sum_{n=0}^{k}\frac{s^{2n}}{n!}
		= \prod_{k=1}^{N-1}\left[ 1 - e^{-s^2}\sum_{n=k+1}^{+\infty}\frac{s^{2n}}{n!}\right]
\end{equation*}
as stated in reference \cite{KhoruzhenkoSommers2015}.\\
\end{appendices}
\newpage
\bibliographystyle{acm}
\bibliography{Library}
\end{document}